\documentclass[11pt]{article}
\usepackage{amsfonts,latexsym,amsmath,amscd,geometry}
\usepackage{amssymb}
\usepackage{latexsym}

\newcommand \nc{\newcommand}
\newtheorem{theorem}{Theorem}[section]
\newtheorem{lemma}[theorem]{Lemma}

\newtheorem{definition}[theorem]{Definition}

\nc{\ba}{\begin{array}}\nc{\ea}{\end{array}}
\nc{\be}{\begin{eqnarray}}\nc{\ee}{\end{eqnarray}}
\nc{\beq}{\begin{equation}}\nc{\eeq}{\end{equation}}
\nc{\bex}{\begin{eqnarray*}}\nc{\eex}{\end{eqnarray*}}
\nc{\btm}{\begin{theorem}} \nc{\etm}{\end{theorem}}
\nc{\blm}{\begin{lemma}} \nc{\elm}{\end{lemma}}
\nc{\R}{\mathbb{R}}  \nc{\ld}{\lambda}
\nc{\va}{\varphi}
\nc{\ve}{\varepsilon}

\def\pf{\noindent{\bf Proof.\quad}}\def\endpf{\hfill$\Box$}

\newcommand \qed {\hfill $\Box$}

\begin{document}
\title{Boundary bubbling analysis of approximate harmonic maps under either weak or strong anchoring  conditions in dimension two}
\author{Tao Huang\footnote{NYU-ECNU Institute of Mathematical Sciences at NYU Shanghai,
3663 Zhongshan Road North, Shanghai, 200062,
China.}, \ \ Changyou Wang\footnote{
Department of Mathematics, Purdue University,
West Lafayette, IN 47907, USA.}}
\maketitle
\begin{abstract}  In this paper, we will study the bubbling phenomena of approximate harmonic maps in dimension two
that have either (i) bounded $L^2$-tension fields under the weak anchoring condition, or (ii) bounded $L\log L \cap M^{1,\delta}$-tension fields under the strong anchoring  condition.
\end{abstract}

\section{Introduction}
\setcounter{equation}{0}
\setcounter{theorem}{0}

The minimization problem of the Landau-De Gennes energy functional for $Q$-tensors under the weak anchoring
boundary conditions has played an important role in the study of nematic liquid crystals (see \cite{De-P, KL, MN, ABG}). It is well-known that
the Landau-De Gennes energy functional for $Q$ reduces to the classical Oseen-Frank energy functional 
for unit vector fields $u$, when $Q$ is restricted to be uniaxial types, i.e., $Q=s(u\otimes u-\frac1{n}\mathbb I_n)$
for a constant scalar order parameter $s$ and a unit vector field $u:\Omega\subset\R^n\to \mathbb S^{n-1}$.
Furthermore, under one constant approximation  the Oseen-Frank energy functional (cf. \cite{HKL}) reduces to the standard
Dirichlet energy functional, whose critical points correspond to harmonic maps. This motivates us to
study the boundary regularity issues of harmonic maps under weak anchoring boundary conditions. 

Let $\Omega\subset \R^n$ be a bounded smooth domain and $N\subset \R^L$ be a compact smooth Riemannian
manifold without boundary.  For a given map $g:\partial\Omega\to N$ and $w>0$, a harmonic map
$u\in H^1(\Omega, N)=\big\{v\in H^1(\Omega, \R^L): \ v(x)\in N \ {\rm{for\ a.e.}}\ x\in \Omega\big\}$, with 
weak anchoring boundary value $g$ and  anchoring strength parameter $w$, if it is a critical point of the modified Dirichlet energy functional:
\begin{equation}\label{MDE}
E(u)=\int_\Omega \frac12|\nabla u|^2+w\int_{\partial\Omega}\frac12 |u-g|^2.
\end{equation}
By direct calculations, we can verify that $u$
solves:
\begin{equation}\label{HM_WA00}
\begin{cases}
\Delta u+A(u)(\nabla u,\nabla u)=0 & \ {\rm{in}}\ \Omega,\\
\frac{\partial u}{\partial\nu}+w\mathbb P(u)(u-g)=0 & \ {\rm{on}}\ \partial\Omega,
\end{cases}
\end{equation}
where $A(\cdot)(\cdot,\cdot)$ denotes the second fundamental form of $N\subset\R^L$,
$\nu$ is the outward unit normal of $\partial\Omega$, and $\mathbb P(y):\R^L\to T_y N$, $y\in N$,
is the orthogonal projection map from $\R^L$ to the tangent space of $N$ at $y$, $T_yN$. 

Recall that a map $u\in H^1(\Omega, N)$ solves the equation \eqref{HM_WA00}, if 
\begin{equation}\label{HM_WA01}
\int_\Omega (-\langle\nabla u,\nabla\phi\rangle+\langle A(u)(\nabla u,\nabla u), \phi\rangle)
+w\int_{\partial\Omega}\langle \mathbb P(u)(u-g),\phi\rangle =0
\end{equation}
holds for all $\phi\in H^1(\Omega,\R^L)\cap L^\infty(\Omega, \R^L)$. Note that when $w=0$, the boundary condition
\eqref{HM_WA00}$_2$ reduces to the zero Neumann boundary condition, which corresponds to
the free anchoring condition; while when $w=\infty$, \eqref{HM_WA00}$_2$ reduces to
the Dirichlet boundary condition $u=g$ on $\partial\Omega$, also called as the strong anchoring condition. 

In a forthcoming paper \cite{LW}, we will extend the interior regularity theorems of
weakly harmonic maps by H\'elein \cite{H} ($n=2$) and stationary harmonic maps
by Bethuel \cite{B} ($n\ge 3$) to the boundary under weak anchoring conditions (see
also the related earlier works on the boundary regularity of harmonic maps
under the Dirichlet boundary condition by Qing \cite{Q} and Wang \cite{W}).
In particular, we will establish the complete boundary regularity of weakly harmonic maps
with weak anchoring condition and a partial boundary regularity for 
stationary harmonic maps in dimensions $n=2$ and $n\ge 3$ respectively.

In this paper, we will mainly be interested in the boundary asymptotic behavior of weakly convergent sequences of (approximate)
harmonic maps with weak anchoring conditions in dimension $n=2$. There have been extensive studies on
the interior asymptotic behavior for sequences of (approximate) harmonic maps, called as the bubble tree convergence,
see for example \cite{Jost, Parker, Q1, Ding-Tian, Wang96, Qing-Tian, Lin-Wang, Lin-Wang1, LT, LS, Li-Zhu, Lamm, Li-Wang1, Li-Wang2, Zhu, 
WWZ, wang16}. However, there are very few works addressing the boundary asymptotic behavior of weakly convergent
(approximate) harmonic maps under various boundary conditions in dimension two. Very recently, there is an interesting
work by Jost-Liu-Zhu \cite{JLZ} that studies the boundary blow-up analysis of approximate harmonic maps under free boundary
conditions. Here we plan to analyze the boundary asymptotic behaviors for such sequences. 
Because of the flexibility of the argument, we can analyze approximate harmonic maps either under weak anchoring conditions, whose tension fields belong to $L^2(\Omega,\R^L)$, or under strong anchoring (or Dirichlet) conditions, whose tension fields
belong to $(L\log L\cap M^{1,a})(\Omega,\R^L)$ for some $0\le a<2$. 

\begin{definition} For a constant $w>0$ and a measurable map $g:\partial\Omega\to N$,
a map $u\in H^1(\Omega, N)$ is called an approximate harmonic map under weak anchoring
condition, with a tension field $\tau\in L^2(\Omega, T_uN)$, if  $u$ is a weak solution of
\begin{equation}\label{AHM_WA00}
\begin{cases}
\Delta u+A(u)(\nabla u,\nabla u)=\tau & \ {\rm{in}}\ \Omega,\\
\frac{\partial u}{\partial\nu}+w\mathbb P(u)(u-g)=0 & \ {\rm{on}}\ \partial\Omega,
\end{cases}
\end{equation}
or, equivalently,
\begin{equation}\label{AHM_WA02}
\int_\Omega (-\langle\nabla u,\nabla\phi\rangle+\langle A(u)(\nabla u,\nabla u), \phi\rangle)
+w\int_{\partial\Omega}\langle \mathbb P(u)(u-g),\phi\rangle =\int_\Omega \langle \tau,\phi\rangle,
\end{equation}
holds for all $\phi\in H^1(\Omega,\R^L)\cap L^\infty(\Omega, \R^L)$.
\end{definition}

For any given function $\tau\in L^2(\Omega, \R^L)$, $w>0$, and a measurable map
$g:\partial\Omega\to N$, it is not hard to check that any critical point $u\in H^1(\Omega,N)$ of 
the energy functional
$$\widetilde{E}(u)=\int_\Omega \frac12|\nabla u|^2+w\int_{\partial\Omega}\frac12|u-g|^2
-\int_\Omega \langle \tau, u\rangle,$$
gives rise an approximate harmonic map under weak anchoring condition, with tension
field $\widetilde{\tau}\equiv\mathbb P(u)(\tau)\in L^2(\Omega, T_uN)$. Furthermore, it is readily seen that 
there always exists at least a minimizer of the energy functional $\widetilde{E}$ over $H^1(\Omega, N)$.

To simplify the analysis, we will assume that approximate harmonic maps under weak anchoring conditions
further belong to $H^{2}(\Omega, N)$, which actually are consequences of the regularity theorems by \cite{LW}.

\begin{theorem}\label{bubbling1}
Assume that $\{u_n\}\subset H^{2}(\Omega, N)$ is a sequence of approximate harmonic maps satisfying
\begin{equation}\label{HM_WA02}
\begin{cases}
\Delta u_n+A(u_n)(\nabla u_n, \nabla u_n)=\tau_n & \ {\rm{in}}\ \Omega,\\
\frac{\partial u_n}{\partial\nu}+w_n\mathbb P(u_n)(u_n-g_n)=0 & \ {\rm{on}}\ \partial\Omega,
\end{cases}
\end{equation}
with $\tau_n\in L^2(\Omega, T_{u_n} N)$, $g_n\in H^{\frac12}(\partial\Omega,N)$, and $w_n>0$.
Assume that there exists $C>0$ such that
\begin{equation}\label{uniform_bound}
\sup_{n\ge 1}\Big\{ \|\nabla u_n\|_{L^2(\Omega)}+\|\tau_n\|_{L^2(\Omega)}+\|g_n\|_{H^{\frac12}(\partial\Omega)}+w_n\Big\}\le C.
\end{equation}
Then there exist a non-negative integer $m$, $u\in H^2(\Omega, N)$, $\tau\in L^2(\Omega, T_uN)$, $w\ge 0$, $g\in H^{\frac12}(\partial\Omega, N)$,
nontrivial harmonic maps $\{\omega_i\}_{i=1}^m\subset C^\infty(\mathbb S^2, N)$, sequences of
points $\{x_n^i\}_{i=1}^m\subset\overline\Omega$, and sequences of scales $\{r_n^i\}_{i=1}^m\subset (0,\infty)$
such that, after passing to a subsequence,
$$u_n\rightharpoonup u \ {\rm{in}}\ H^1(\Omega),\ \tau_n\rightharpoonup \tau \ {\rm{in}}\ L^2(\Omega),
\ w_n\rightarrow w, 
\ g_n\rightharpoonup g \ {\rm{in}}\ H^{\frac12}(\partial\Omega).
$$
Moreover, the following statements hold:
\begin{itemize}
\item[(i)] $u\in H^2(\Omega, N)$ is an approximate harmonic map with tension field $\tau$ with weak anchoring
condition $g$ and anchoring strength parameter $w$, i.e., satisfying 
\begin{equation}\label{HM_WA03}
\begin{cases}
\Delta u+A(u)(\nabla u, \nabla u)=\tau & \ {\rm{in}}\ \Omega,\\
\frac{\partial u}{\partial\nu}+w\mathbb P(u)(u-g)=0 & \ {\rm{on}}\ \partial\Omega,
\end{cases}
\end{equation}

\item[(ii)] For any pair $1\le i<j\le m$, 
\begin{equation}\label{bubble-scales}
\lim_{n\rightarrow\infty}\Big\{\frac{r_n^i}{r_n^j}+\frac{r_n^j}{r_n^i}+\frac{|x_n^i-x_n^j|}{r_n^i+r_n^j}\Big\}=\infty.
\end{equation}
\item[(iii)] {\rm{(energy identity)}} 
\begin{equation}\label{energy_id}
\lim_{n\rightarrow\infty}\|\nabla u_n\|_{L^2(\Omega)}^2
=\|\nabla u\|_{L^2(\Omega)}^2+\sum_{i=1}^m \|\nabla \omega_i\|_{L^2(\Omega)}^2.
\end{equation}
\item[(iv)] {\rm{(oscillation convergence)}}
\begin{equation}\label{oscillation_conv}
\lim_{n\rightarrow\infty}\Big\|u_n-u-\sum_{i=1}^m \big(\omega_i(\frac{\cdot-x_n^i}{r_n^i})-\omega_i(\infty)\big)\Big\|_{L^\infty(\Omega)}=0.
\end{equation}
\end{itemize}
\end{theorem}

A few remarks are in order:

\smallskip
\noindent
(1) It follows from Theorem \ref{bubbling1} that  when $\tau_n=0$, the bubble tree convergence holds for harmonic maps
under the weak anchoring condition in dimension two.

\smallskip
\noindent
(2) In a forthcoming paper \cite{Huang-Wang}, we will establish the existence of a global weak solution
to the heat flow of harmonic maps under weak anchoring conditions in dimension two, extending the works
by Struwe \cite{Struwe} and Chang \cite{chang89}, and discuss the application of Theorem \ref{bubbling1}
to the heat flow of harmonic maps under weak anchoring condition in dimension two.

\smallskip
\noindent
(3) It remains to be an interesting question to ask whether Theorem \ref{bubbling1} holds if we relax the integrability of
$\tau_n$ to the class that $\tau_n\in L\log L\cap M^{1,a}$, for some $0\le a<2$, are bounded. Here $M^{1,a}$ denotes the Morrey space ($1,a$). The interior case of Theorem \ref{bubbling1} does hold when $\tau_n$ is bounded
in $L\log L\cap \mathcal M^{1,a}$ for some $0\le a<2$ by a recent work by the author \cite{wang16} (see also
a related work \cite{WWZ}).

To carry out the boundary blowing up analysis for approximate harmonic maps $u_n$ with $L^2$-tensions under weak anchoring
conditions, we first need to establish a boundary H\"older continuity estimate under the small energy assumption.
This involves several new observations: (1) using the star-shape property of the half ball $B_r(x)\cap\Omega$, with $x\in\partial\Omega$
and small $r>0$, we can apply a Pohozaev type argument, see Lemma 2.2,  to control the oscillation of $u_n$ on $B_r(x)\cap\Omega$; (2) the Courant-Lebesgue Lemma can be used to control the oscillation of $u_n$ on $\partial B_r(x)\cap\Omega$;
and these two ingredients, combined with the interior H\"older continuity estimate, can yield the boundary H\"older estimate
of $u_n$, see Theorem 2.3. Second, it follows from simple scaling arguments that the blowing up limit of the weak anchoring condition is zero Neumann condition, and hence any bubble at a boundary concentration point is also a harmonic $\mathbb S^2$. The most difficult step is to show the vanishing of energy and oscillation in a boundary neck region. This involves to establish that the energy of $u_n$ on any dyadic boundary annual $\big(B_{2r}(x)\setminus B_r(x)\big)\cap\Omega$, $x\in\partial\Omega$, within the boundary neck region decays like $r^\alpha$ for $\alpha\in(0,1)$, see Lemma 3.1. A crucial step here is to control the radial
energy of $u_n$ on $\partial B_r(x)\cap\Omega$, $x\in\partial\Omega$, by the tangential energy of $u_n$  on $\partial B_r(x)\cap\Omega$ along with the bulk energy of $u_n$ and $L^2$-energy of tension fields $\tau_n$ on $B_r(x)\cap\Omega$ and $H^\frac12$-norm of $g_n$ on $\partial\Omega\cap B_r(x)$, see Lemma 2.4.

\smallskip
Since Theorem \ref{bubbling1} requires that the anchoring strength parameters $w_n$ are uniformly bounded, it does not 
apply to the case of strong anchoring condition or the Dirichlet boundary condition. However, the global bubble tree convergence remains to be true for approximate harmonic maps under strong anchoring conditions. 
To state the result, we recall a few notations. The space $L\log L(\Omega)$ is defined by 
$$
L\log L(\Omega)
=\left\{
f\in L^1(\Omega):\big\|f\big\|_{L\log L(\Omega)}=\int_{\Omega}|f|\log(2+|f|)<+\infty
\right\},
$$
and the Morrey space $\mathcal M^{p,a}(\Omega)$, for $1\leq p<+\infty$ and $0\leq a\leq 2$, is defined by 
$$
\mathcal M^{p,a}(\Omega)
=\left\{
f\in L^{p}_{\mbox{loc}}(\Omega):\|f\|^p_{\mathcal M^{p,a}(\Omega)}=\sup\limits_{B_r(x)\subset\Omega}r^{a-2}\int_{B_r(x)}|f|^p<\infty
\right\}.
$$ 
The strong anchoring boundary data $h_n:\partial\Omega\rightarrow N$ is assumed to
satisfy the following two assumptions:
\begin{itemize}
\item[(A1)] $\{h_n\}$ is uniformly continuous on $\partial\Omega$: for any $\epsilon>0$, there exists  $\delta>0$ such that 
$\displaystyle\sup_{n\ge 1}|h_n(x)-h_n(y)|\leq \epsilon$, whenever $x,y\in\partial\Omega$ satisfies $|x-y|\leq \delta$.

\item[(A2)] $\{h_n\}$ is equi-integrable in $H^{\frac12}(\partial\Omega)$ in the sense that
$$
\lim_{E\subset\partial\Omega, \ H^1(E)\rightarrow 0}\sup\limits_{n}\big\|h_n\big\|_{H^{\frac12}(E)}=0.
$$
\end{itemize}
The global bubble tree convergence result for approximate harmonic maps under the Dirichlet boundary condition
can be stated as follows.
\btm\label{mainth1}
Assume $\{h_n\}\subset C^0(\partial\Omega,N)\cap H^{\frac12}(\partial\Omega,N)$ satisfies the assumptions (A1) and (A2). Let $\{u_n\}\subset H^1(\Omega, N)$ be a sequence of approximate harmonic maps under the Dirichlet boundary condition: 
\beq\label{nahm}
\left\{
\ba{rll}
\Delta u_n+A(u_n)(\nabla u_n,\nabla u_n)=&\tau_n,\quad &\mbox{in } \Omega,\\
u_n=&h_n,\quad &\mbox{on }\partial\Omega,
\ea
\right.
\eeq
satisfying
\beq\label{theq1}
\|\nabla u_n\|_{L^2(\Omega)}+\|\tau_n\|_{L\log L(\Omega)}
+\|\tau_n\|_{\mathcal M^{1,a}(\Omega)}\leq C<+\infty
\eeq
for some fixed $0\leq a<2$. Then, after taking a subsequence, we have that
$$
u_n\rightharpoonup u \mbox{ in }H^1(\Omega),\quad
\tau_n\rightharpoonup \tau \mbox{ in }L^1(\Omega),\quad 
h_n\rightarrow h \mbox{ in }C^0(\partial\Omega). 
$$
Moreover, the following statements hold:
\begin{enumerate}
\item The limit function $u\in H^1(\Omega, N)\cap C^0(\overline\Omega,N)$ is an approximate harmonic map,
with tension field $\tau\in L\log L(\Omega)\cap \mathcal M^{1,a}(\Omega)$ and strong anchoring condition $h
\in C^0(\partial\Omega,N)\cap H^{\frac12}(\partial\Omega,N)$:
\beq\notag
\left\{
\ba{rll}
\Delta u+A(u)(\nabla u,\nabla u)=&\tau,\quad &\mbox{in } \Omega,\\
u=&h,\quad &\mbox{on }\partial\Omega.
\ea
\right.
\eeq
\item There exist  a nonnegative integer $m$, sequences of points $\{x_n^i\}_{i=1}^m\subset\bar\Omega$, 
sequences of scales $\{r_n^i\}_{i=1}^m\subset\R_+$, with $r_n^i\downarrow 0$ for any $1\leq i\leq m$, 
and nontrivial harmonic maps $\{\omega_i\}_{i=1}^{m}\in C^\infty(\mathbb S^2, N)$ such that
\beq\label{theq2}
\lim_{n\rightarrow\infty}\big\{\frac{r_n^i}{r_n^j}+\frac{r_n^j}{r_n^i}+\frac{|x_n^i-x_n^j|}{r_n^i+r_n^j}\big\}=\infty, 
\ \forall\ 1\leq i< j\leq m,
\eeq
\beq\label{theq3}
\lim\limits_{n\rightarrow \infty}\|\nabla u_n\|^2_{L^2(\Omega)}=\|\nabla u\|^2_{L^2(\Omega)}
+\sum\limits_{i=1}^m\int_{\mathbb S^2}|\nabla \omega_i|^2,
\eeq
and 
\beq\label{theq4}
\lim_{n\rightarrow\infty}\big\|u_n-u-\sum\limits_{j=1}^l\omega_n^j\big\|_{L^\infty(\Omega)}=0,
\ {\rm{with}}\ \omega_n^j(\cdot)=\omega_j\big(\frac{\cdot-x_n^j}{r_n^j}\big)-\omega_j(\infty)\ (1\le j\le m).
\eeq
\end{enumerate}
\etm

The ideas to prove Theorem \ref{mainth1} are similar yet much simpler than that of Theorem \ref{bubbling1}.
First, we can show a uniform boundary H\"older continuity estimate of $u_n$ under a small energy assumption, see Theorem 6.2. Second,
observe that any bubble at a boundary concentration point is a harmonic $\mathbb S^2$, because the blowing up limit of
$g_n$ is constant. Third, when the boundary data $g_n$ is uniformly continuous, one can rather easily show there is no oscillation accumulation in any boundary neck region, which can then be used to show no energy concentration in any boundary neck region, see Lemma 7.1.

The paper is organized as follows. In section 2,  we will establish the boundary H\"older continuity estimate of approximate harmonic maps with weak anchoring conditions, under the small energy condition; and a boundary Rellich's inequality controlling
radial energy by tangential energy. In section 3, we will estimate both energy and oscillation of approximate harmonic maps
with weak anchoring conditions in any boundary neck region. In section 4, we will prove a $H^2$-type removability
of an isolated boundary singularity. In section 5, we give a proof of Theorem \ref{bubbling1}. The section 6
is devoted to  the boundary H\"older continuity estimate of approximate harmonic maps with Dirichlet conditions, under
the small energy condition. The section 7 is devoted
to the proof of Theorem \ref{mainth1}.

\section{Some lemmas for approximate harmonic maps with weak anchoring conditions}
\setcounter{equation}{0}
\setcounter{theorem}{0}

In this section, we will establish some crucial lemmas that are needed in the proof of Theorem \ref{bubbling1}.
For $x_0\in \Omega$ and $r>0$, we denote by $B_r(x_0)$ the ball in $\R^2$ with center $x_0$ and radius $r$, and
$$B_r^+(x_0)=B_r(x_0)\cap\Omega, \ S_r^+(x_0)=\partial B_r(x_0)\cap\Omega,
\ T_r(x_0)=B_r(x_0)\cap\partial\Omega.$$

We begin with an interior estimate for approximate harmonic maps.
\begin{lemma}\label{interior}  There exists $\epsilon_0>0$ such that if $u\in H^2(\Omega, N)$ 
is an approximate harmonic map with tension field $\tau\in L^2(\Omega, T_u N)$, and satisfies, for some ball $B_r(x_0)\subset\Omega$,
$$\int_{B_r(x_0)}|\nabla u|^2\le\epsilon_0^2,$$
then
\begin{equation}\label{H2_bound}
r^2\|\nabla^2 u\|_{L^2(B_{\frac{r}2}(x_0))}^2
\le C\big(r^2\|\tau\|_{L^2(B_r(x_0))}^2+\|\nabla u\|_{L^2(B_r(x_0))}^2\big).
\end{equation}
In particular, $u\in C^\gamma(B_{\frac{r}2}(x_0))$ for all $\gamma\in (0,1)$, and
\begin{equation}
\label{holder_estimate}
{\rm{osc}}_{B_s(x_0)} u
\le C(\gamma)(\frac{s}{r})^\gamma\big(r\|\tau\|_{L^2(B_r(x_0))}+\|\nabla u\|_{L^2(B_r(x_0))}\big), \ \forall \ 0<s\le \frac{r}2.
\end{equation}
\end{lemma}
\pf \eqref{H2_bound} follows from suitable modifications of the argument by Sacks-Uhlenbeck \cite{SaU}, see Ding-Tian \cite{Ding-Tian} for the details. \eqref{holder_estimate} follows from \eqref{H2_bound} and Sobolev's embedding theorem. \hfill\qed

\medskip
To deal with the weak anchoring  condition \eqref{HM_WA00}$_2$, we need to establish boundary estimates
analogous to that of Lemma \ref{interior}. To do it, we first control the tangential energy of $u$ on $\partial\Omega$ locally, which invokes a local nonlinear version of Rellich's type argument.

\begin{lemma}\label{rellich1} Assume that $u\in H^2(\Omega, N)$ solves 
the equation \eqref{HM_WA00}, with $\tau\in L^2(\Omega, T_uN)$, $w>0$, 
and $g\in H^{\frac12}(\partial\Omega, N)$. Then there
exists $r_0=r_0(\partial\Omega)>0$ such that for any $x_0\in\partial\Omega$ and 
$0<r\le r_0$, it holds that
\begin{equation}\label{rellich2} 
r\int_{T_r(x_0)}|\nabla_T u|^2\le
C\big\{\int_{B_{2r}^+(x_0)}(|\nabla u|^2+r^2|\tau|^2)+w^2 r\int_{T_{2r}(x_0)}|u-g|^2\big\}.
\end{equation}
Here $\nabla_T$ denotes the tangential derivative on $\partial\Omega$. In particular,
we have that
\begin{equation}\label{bdry_oscillation}
\big({\rm{osc}}_{T_r(x_0)} u\big)^2
\le C\big\{\int_{B_{2r}^+(x_0)}(|\nabla u|^2+r^2|\tau|^2)+w^2 r\int_{T_{2r}(x_0)}|u-g|^2\big\}.
\end{equation}
\end{lemma}
\pf Since $\partial\Omega$ is smooth, it is well-known that there exists $r_0=r_0(\partial\Omega)>0$ such that
for any $x_0\in\partial\Omega$ and $0<r\le r_0$, $B_r^+(x_0)$ is star-shaped with a center $a\in B_r^+(x_0)$
in the sense that there exists a universal positive constant $c_0$ such that
\begin{equation} \label{star-shape}
(x-a)\cdot\nu\ge c_0 r, \ \forall x\in \partial B_r^+(x_0),
\end{equation}
where $\nu$ is the outward unit normal of $\partial B_r^+(x_0)$.

For simplicity, we may further assume, by Fubini's theorem,  that 
\begin{equation}\label{fubini}
r\int_{S_r^+(x_0)}|\nabla u|^2\le 8\int_{B_{2r}^+(x_0)}|\nabla u|^2.
\end{equation}
Multiplying \eqref{HM_WA00}$_1$ by $(x-a)\cdot\nabla u$ and integrating the resulting equation over
$B_r^+(x_0)$, we obtain
\begin{eqnarray*}
&&\int_{B_r^+(x_0)}\langle\tau, (x-a)\cdot\nabla u\rangle\\
&&=\int_{B_r^+(x_0)}\langle\Delta u, (x-a)\cdot\nabla u\rangle\\
&&=\int_{B_r^+(x_0)}{\rm{div}}\langle\nabla u, (x-a)\cdot\nabla u\rangle-\int_{B_r^+(x_0)}|\nabla u|^2
-\frac12\int_{B_r^+(x_0)}(x-a)\cdot\nabla(|\nabla u|^2)\\
&&=\int_{\partial B_r^+(x_0)}\langle\frac{\partial u}{\partial\nu}, (x-a)\cdot\nabla u\rangle
-\frac12\int_{\partial B_r^+(x_0)}(x-a)\cdot\nu |\nabla u|^2.
\end{eqnarray*}
This, with the help of \eqref{star-shape}, Young's inequality, H\"older's inequality, and \eqref{fubini}, 
implies that
\begin{eqnarray*}
&&r\int_{\partial B_r^+(x_0)}|\nabla u|^2\\
&&\le Cr\int_{\partial B_r^+(x_0)}|\frac{\partial u}{\partial\nu}|^2+C\int_{B_r^+(x_0)}|\tau||x-a||\nabla u|\\
&&\le Cr\int_{S_r^+(x_0)}|\frac{\partial u}{\partial\nu}|^2+Cw^2r\int_{T_r(x_0)}|u-g|^2
+C\int_{B_r^+(x_0)}|\nabla u|^2+Cr^2\int_{B_r^+(x_0)}|\tau|^2\\
&&\le C\int_{B_{2r}^+(x_0)}|\nabla u|^2+Cr^2\int_{B_r^+(x_0)}|\tau|^2+Cw^2r\int_{T_r(x_0)}|u-g|^2.
\end{eqnarray*}
Since $\int_{T_r(x_0)}|\nabla_T u|^2\le \int_{\partial B_r^+(x_0)}|\nabla u|^2$, this clearly yields \eqref{rellich2}. 
\eqref{bdry_oscillation} follows from \eqref{rellich2} and the following inequality
$$|u(x)-u(y)|\le \int_{T_r(x_0)}|\nabla_T u|
\le Cr^\frac12\big(\int_{T_r(x_0)}|\nabla_T u|^2\big)^\frac12, \ \forall x, y\in T_r(x_0).$$
The proof is complete.
\qed

\medskip
With the help of Lemma \ref{interior} and Lemma \ref{rellich1}, we can prove a local boundary oscillation estimate of
approximate harmonic maps under the weak anchoring boundary condition. More precisely,
we have
\begin{theorem} \label{bdry_max1} There exist $\epsilon_0>0$ and  $r_0=r_0(\partial\Omega)>0$
such that if $u\in H^2(\Omega, N)$ solves the equation \eqref{HM_WA00},
with  $\tau\in L^2(\Omega, T_uN)$, $w>0$, and $g\in H^{\frac12}(\partial\Omega, N)$, 
and satisfies, for some $x_0\in\partial\Omega$ and $0\le r<r_0$, 
$$\int_{B_{2r}^+(x_0)}|\nabla u|^2\le\epsilon_0^2,$$
then it holds that
\begin{equation}\label{bdry_max2}
\big({\rm{osc}}_{B_{\frac{r}2}^+(x_0)} u\big)^2
\le C\Big\{\int_{B_{2r}^+(x_0)}(|\nabla u|^2+r^2|\tau|^2)+w^2 r\int_{T_{2r}(x_0)}|u-g|^2\Big\},
\end{equation}
and
\beq\label{H2-estimate}
r^2\big\|\nabla^2u\big\|^2_{L^2(B_{\frac{r}2}^+(x_0))}\le
C\Big\{\int_{B_{r}^+(x_0)}((1+w^2)|\nabla u|^2+r^2|\tau|^2)+w^2r^2\|g\|^2_{H^\frac12(T_{r}(x_0))}\Big\}.
\eeq
\end{theorem}
\pf By Fububi's theorem, we may assume, for simplicity, that 
\begin{equation}\label{fubini1}
r\int_{S_r^+(x_0)}|\nabla u|^2\le 8\int_{B_{2r}^+(x_0)}|\nabla u|^2.
\end{equation}
This, combined with the Sobolev's embedding theorem, implies that $u\in C^\frac12(S_r^+(x_0))$ and
\begin{equation}\label{fubini_osc}
{\rm{osc}}_{S_r^+(x_0)} u\le C\big(r\int_{S_r^+(x_0)}|\nabla u|^2\big)^\frac12
\le C\big(\int_{B_{2r}^+(x_0)}|\nabla u|^2\big)^\frac12.
\end{equation}
This, combined with \eqref{bdry_oscillation}, implies that
\begin{equation}\label{bdry_oscillation1}
\big({\rm{osc}}_{\partial B_r^+(x_0)} u\big)^2
\le C\int_{B_{2r}^+(x_0)}(|\nabla u|^2+r^2|\tau|^2)+Cw^2 r\int_{T_r(x_0)}|u-g|^2.
\end{equation}
Now we want to show that there exists $C>0$ such that for any $P\in N$, 
\begin{equation}\label{local_max}
\sup_{B_{\frac{r}2}^+(x_0)}|u(x)-P|\le \sup_{\partial B_r^+(x_0)}|u(x)-P|+C
\big\{\int_{B_{2r}^+(x_0)}(|\nabla u|^2+r^2|\tau|^2)\big\}^\frac12.
\end{equation}
It is readily seen that \eqref{bdry_max2} follows directly from \eqref{bdry_oscillation1}
and \eqref{local_max}, since $P\in N$ can be an arbitrary point. 

To prove \eqref{local_max}, set 
$$L=\sup_{B_{\frac{r}2}^+(x_0)}|u(x)-P|<+\infty,$$
and choose $x_1\in B_{\frac{r}2}^+(x_0)$ so that
$$|u(x_1)-P|\ge \frac{L}2.$$
We may further assume that
 $$L\ge 64 \big\{\int_{B_{2r}^+(x_0)}(|\nabla u|^2+r^2|\tau|^2)\big\}^\frac12.$$
 Set $d_1={\rm{dist}}(x_1, \partial B_{\frac{r}2}^+(x_0))>0$. Then $B_{d_1}(x_1)\subset \Omega\cap B_{2r}^+(x_0)$,
 and Lemma \ref{interior} implies that
 for any $0<\theta<1$, 
 $${\rm{osc}}_{B_{\theta d_1}(x_1)} u\le C\theta^\frac12 
 \big\{\int_{B_{d_1}(x_1)}(|\nabla u|^2+r^2|\tau|^2)\big\}^\frac12\le C\theta^\frac12 L.$$
 Choosing $\theta_0=\frac{1}{16C^2}\in (0,1)$, we obtain that
 \begin{equation}\label{lower_bound}
 \inf_{B_{\theta_0 d_1}(x_1)}|u(x)-P|\ge \frac{L}{4}.
 \end{equation}
 Let $\widehat{u}: B_{2r}(x_0)\to \R^L$ be  such that
 $\widehat u=u$ in $B_{2r}^+(x_0)$, and 
 $$\int_{B_{2r}(x_0)}|\nabla\widehat{u}|^2 \le 4\int_{B_{2r}^+(x_0)}|\nabla u|^2.$$
 By Fubini's theorem, there exists $r_1\in (d_1, 2d_1)$ such that
 $$r_1\int_{\partial B_{r_1}(x_1)}|\nabla\widehat u|^2\le 4\int_{B_{2r}(x_0)}|\nabla \widehat u|^2
 \le 16\int_{B_{2r}^+(x_0)}|\nabla u|^2,
 $$
 and hence, by Sobolev's embedding theorem, 
 $${\rm{osc}}_{\partial B_{r_1}(x_1)} {\widehat u}\le \big(r_1\int_{\partial B_{r_1}(x_1)}|\nabla\widehat u|^2\big)^\frac12
 \le (16\int_{B_{2r}^+(x_0)}|\nabla u|^2)^\frac12\le \frac{L}{16}.$$
 Since $\partial B_{r_1}(x_1)\cap \partial B_{r}^+(x_0)\not=\emptyset$, this yields
 \begin{equation}\label{upper_bound}
 \sup_{\partial B_{r_1}(x_1)}|\widehat u(x)-P|\le \sup_{\partial B_r^+(x_0)}|u-P|
 +\frac{L}{16}.
 \end{equation}
 It follows from \eqref{lower_bound} and \eqref{upper_bound} that 
 \begin{eqnarray*}
&& \frac{L}4\le \inf_{x\in\partial B_{\theta_0 d_1}(x_1)}|u(x)-P|\\
&&\le \inf_{x\in\partial B_{\theta_0 d_1}(x_1), y\in \partial B_{r_1}(x_1)}|\widehat{u}(x)-\widehat{u}(y)|
+\sup_{y\in \partial B_{r_1}(x_1)}|\widehat{u}(y)-P|\\
&&\le \inf_{x\in\partial B_{\theta_0 d_1}(x_1), y\in \partial B_{r_1}(x_1)}|\widehat{u}(x)-\widehat{u}(y)|
+\sup_{\partial B_r^+(x_0)}|u-P|
 +\frac{L}{16}.
 \end{eqnarray*}
 This implies that
 \begin{equation}\label{upper_bound1}
 \frac{L}{8}\le \inf_{x\in\partial B_{\theta_0 d_1}(x_1), y\in \partial B_{r_1}(x_1)}|\widehat{u}(x)-\widehat{u}(y)|
+\sup_{\partial B_r^+(x_0)}|u-P|.
\end{equation}
Observe that, by using polar coordinates,  it holds 
\begin{eqnarray*}
&&\inf_{x\in\partial B_{\theta_0 d_1}(x_1), y\in \partial B_{r_1}(x_1)}|\widehat{u}(x)-\widehat{u}(y)|\\
&&\le\frac{1}{4\pi^2} \int_{\mathbb S^1}|\widehat{u}(x_1+\theta_0 d_1\theta)-\widehat{u}(x_1+r_1\theta)|\,d\theta\\
&&\le\frac{1}{4\pi^2} \int_{\mathbb S^1}\int_{\theta_0 d_1}^{2d_1}\big|\frac{\partial \widehat{u}}{\partial r}\big|\,drd\theta\\
&&\le \frac{1}{4\pi^2}\big(\int_{\theta_0 d_1}^{2d_1}\frac{1}{r}\,dr\big)^\frac12\int_{\mathbb S^1}\big(\int_{\theta_0 d_1}^{2d_1}\big|\frac{\partial \widehat{u}}{\partial r}\big|^2r\,dr)^\frac12\,d\theta\\
&&\le \frac{1}{4\pi^2}\big(\ln \frac{2}{\theta_0}\big)^\frac12(2\pi)^\frac12 \big(\int_{B_{2r}(x_0)}|\nabla\widehat u|^2\big)^\frac12\\
&&\le C\big(\int_{B_{2r}(x_0)}|\nabla\widehat u|^2\big)^\frac12
\le C\big(\int_{B_{2r}^+(x_0)}|\nabla u|^2\big)^\frac12.
\end{eqnarray*} 
Putting this estimate into \eqref{upper_bound1}, we obtain that
 \begin{equation}\label{upper_bound2}
 \frac{L}{8}\le C\big(\int_{B_{2r}^+(x_0)}|\nabla u|^2\big)^\frac12
+\sup_{\partial B_r^+(x_0)}|u-P|.
\end{equation}
It is readily seen that \eqref{upper_bound2} implies \eqref{local_max}.

To show \eqref{H2-estimate}, choose $r_1\in (\frac{r}2, r)$ such that 
\beq\label{good_r}
r_1\int_{S_{r_1}^+(x_0)}|\nabla ^2 u|^2\le \int_{B_r^+(x_0)}|\nabla^2 u|^2.
\eeq
Let $v\in H^2(B_{r_1}^+(x_0),\R^L)$ solve
\beq\label{auxi-eqn}
\begin{cases}
\Delta v=0 & \ {\rm{in}}\ B_{r_1}^+(x_0),\\
v=u& \ {\rm{on}}\ S_{r_1}^+(x_0),\\
\frac{\partial v}{\partial\nu}=-w(u-g) & \ {\rm{on}}\ T_{r_1}(x_0).
\end{cases}
\eeq
Then by the standard theory for Laplace equations we have that $v\in H^2(B_{r_1}^+(x_0))$, and
\begin{eqnarray}\label{H2-est1}
r_1^2\big\|\nabla^2 v\big\|^2_{L^2(B_{r_1}^+(x_0))}
&\le& C(r_1^3\int_{S_{r_1}^+(x_0)}|\nabla^2 u|^2+w^2r_1^2\|u-g\|^2_{H^{\frac12}(T_{r_1}(x_0))}\big)\nonumber\\
&\le& C((1+w^2)r^2\int_{B_{r}^+(x_0)}|\nabla^2 u|^2
+w^2r^2\|g\|^2_{H^{\frac12}(T_{r}(x_0))}\big).
\end{eqnarray}
It is readily seen that $u-v$ solves
\beq\label{auxi-eqn1}
\begin{cases}
-\Delta (u-v)=A(u)(\nabla u,\nabla u)+\tau & \ {\rm{in}}\ B_{r_1}^+(x_0),\\
u-v=0& \ {\rm{on}}\ S_{r_1}^+(x_0),\\
\frac{\partial }{\partial\nu}(u-v)=0 & \ {\rm{on}}\ T_{r_1}(x_0).
\end{cases}
\eeq
Hence, by the $H^2$-theory and Ladyzhenskaya's inequality, we obtain that
\begin{eqnarray}\label{H2-est2}
&&r_1^2\big\|\nabla^2(u-v)\big\|^2_{L^2(B_{r_1}^+(x_0))}
\le Cr_1^2\int_{B_{r_1}^+(x_0)}(|\nabla u|^4+|\tau|^2)\nonumber\\
&&\le C\int_{B_{r_1}^+(x_0)}|\nabla u|^2
\big(\int_{B_{r_1}^+(x_0)}|\nabla u|^2+r_1^2\int_{B_{r_1}^+(x_0)}|\nabla^2 u|^2\big)
+Cr_1^2\int_{B_{r_1}^+(x_0)}|\tau|^2\nonumber\\
&&\le C\epsilon_0^2\big(\int_{B_{r_1}^+(x_0)}|\nabla u|^2+r_1^2\int_{B_{r_1}^+(x_0)}|\nabla^2 u|^2\big)
+Cr_1^2\int_{B_{r_1}^+(x_0)}|\tau|^2
\end{eqnarray}
Adding \eqref{H2-est1} with \eqref{H2-est2} and choosing sufficiently small $\epsilon_0>0$, we then obtain
$$
r_1^2\big\|\nabla^2u\big\|^2_{L^2(B_{r_1}^+(x_0))}
\le C(\epsilon_0^2+w^2)\int_{B_{r}^+(x_0)}|\nabla u|^2+Cr^2\int_{B_{r}^+(x_0)}|\tau|^2
+w^2r^2\|g\|^2_{H^{\frac12}(T_{r}(x_0))}.
$$
This clearly implies \eqref{H2-estimate}. Hence the proof is complete.\qed

\medskip
Finally, we need to control the radial energy of $u$ on $S_r^+(x_0)$ 
by the tangential energy on $S_r^+(x_0)$ for $x_0\in\partial\Omega$ and $r>0$. 
More precisely, we have
\begin{lemma}\label{radial_tangent1} There exist $r_0=r_0(\partial\Omega)>0$, and
$C_0>0$ depending on $w, g,\partial\Omega$ such that
if $u\in H^2(\Omega, N)$ solves the equation \eqref{HM_WA00}, with $\tau\in L^2(\Omega, T_uN)$, $w>0$,
and $g\in H^{\frac12}(\partial\Omega, N)$. 
Then for any $x_0\in\partial\Omega$ and $0<r\le r_0$, it holds
\begin{equation}\label{radial_tangent2}
r\int_{S_r^+(x_0)}|\frac{\partial u}{\partial\nu}|^2
\le Cr\big(1+\big\|g\big\|^2_{H^{\frac12}(\partial\Omega)}+\int_{\Omega}|\nabla u|^2+\int_{S_r^+(x_0)}|\nabla_T u|^2
+\int_{B_r^+(x_0)}|\tau||\nabla u|\big).
\end{equation}
\end{lemma}

\pf For any smooth vector field $X\in C^\infty(\overline\Omega,\R^2)$, 
multiplying \eqref{HM_WA00} by $X\cdot\nabla u$ and integrating the resulting equation
over $B_r^+(x_0)$, we obtain that
\begin{eqnarray}\label{pohozaev}
&&\int_{B_r^+(x_0)}\langle\tau, X\cdot\nabla u\rangle\nonumber\\
&&=\int_{B_r^+(x_0)}\langle\Delta u, X\cdot\nabla u\rangle\nonumber\\
&&=\int_{B_r^+(x_0)}{\rm{div}}\langle X\cdot\nabla u, \nabla u\rangle
-\int_{B_r^+(x_0)}X\cdot\nabla (\frac{|\nabla u|^2}2)-\int_{B_r^+(x_0)} \nabla u\otimes\nabla u: \nabla X\nonumber\\
&&=\int_{\partial B_r^+(x_0)}\langle X\cdot\nabla u, \frac{\partial u}{\partial\nu}\rangle
-\int_{\partial B_r^+(x_0)}X\cdot\nu \frac{|\nabla u|^2}2\nonumber\\
&&\ \ \ +\frac12 \int_{B_r^+(x_0)}\big(|\nabla u|^2{\rm{div}} X-2\nabla u\otimes\nabla u: \nabla X\big).
\end{eqnarray}
Since $\partial\Omega$ is smooth, there exist $r_0=r_0(\partial\Omega)>0$
and $C_0=C_0(\partial\Omega)>0$ and a vector field $X\in C^\infty(B_{r_0}(x_0), \R^2)$
such that the following properties hold:
\begin{equation}\label{nice-vf}
\begin{cases}
X\cdot\nu =0 \ {\rm{on}}\ T_{r_0}(x_0),\\
|X(x)-(x-x_0)|\le C_0|x-x_0|^2, \  \forall x\in B_{r_0}(x_0),\\
|\nabla X(x)-\mathbb I_2|\le 2C_0|x-x_0|, \ \forall x\in B_{r_0}(x_0).
\end{cases}
\end{equation}
For $0<r<r_0$, substituting this $X$ into \eqref{pohozaev} and applying
the boundary condition \eqref{HM_WA00}$_2$, we obtain that
\begin{eqnarray*}
&&(1-Cr)r\int_{S_r^+(x_0)}|\frac{\partial u}{\partial \nu}|^2\\
&&\le (1+Cr)r\int_{S_r^+(x_0)}\frac12|\nabla u|^2
+Cr\int_{B_r^+(x_0)}|\nabla u|^2+(1+C)r\int_{B_r^+(x_0)}|\tau||\nabla u|\\
&&\ +w\int_{T_r(x_0)} \langle X\cdot\nabla_T u, u-g\rangle.
\end{eqnarray*}
We can estimate 
\begin{eqnarray*}
\big|\int_{T_r(x_0)} \langle X\cdot\nabla_T g, u-g\rangle\big|
&\le& r\|\nabla_T g\|_{H^{-\frac12}(\partial\Omega)}\|u-g\|_{H^\frac12(\partial\Omega)}\\
&\le& Cr\|g\|_{H^\frac12(\partial\Omega)}\|u-g\|_{H^\frac12(\partial\Omega)}\\
&\le& Cr\big(\|g\|^2_{H^\frac12(\partial\Omega)}+\|u\|^2_{H^\frac12(\partial\Omega)}\big)\\
&\le&Cr\big(1+\|g\|^2_{H^\frac12(\partial\Omega)}+\|\nabla u\|^2_{L^2(\Omega)}\big),
\end{eqnarray*}
and hence
\begin{eqnarray*}
\big|\int_{T_r(x_0)} \langle X\cdot\nabla_T u, u-g\rangle\big|
&=&\big|\int_{T_r(x_0)} \langle X\cdot\nabla_T (u-g), u-g\rangle
+ \int_{T_r(x_0)} \langle X\cdot\nabla_T g, u-g\rangle\big| \\
&=&\big|\int_{T_r(x_0)} X\cdot\nabla_T(\frac{|u-g|^2}2)
+ \int_{T_r(x_0)} \langle X\cdot\nabla_T g, u-g\rangle\big| \\
&=&\big|-\frac12\int_{T_r(x_0)}{\rm{div}}_{T_r(x_0)}(X)|u-g|^2+\frac12\int_{\partial T_r(x_0)}X\cdot\nu_{\partial T_r(x_0)} |u-g|^2\\
&&\  +\int_{T_r(x_0)} \langle X\cdot\nabla_T g, u-g\rangle\big|\\
&\le& C\big(1+\|g\|^2_{H^{\frac12}(\partial\Omega)}+\|\nabla u\|^2_{L^2(\Omega)}\big)r,
\end{eqnarray*}
where ${\rm{div}}_{T_r(x_0)}(X)$ denotes the divergence of $X$ with respect to $T_r(x_0)$ and
$\nu_{\partial T_r(x_0)}$ denotes the outward unit normal of $\partial T_r(x_0)$. 

Therefore, by choosing sufficiently small $r_0>0$, we conclude that
\begin{equation}\label{radial_tangent_estimate}
r\int_{S_r^+(x_0)}|\frac{\partial u}{\partial\nu}|^2\le Cr\big(1+\big\|g\big\|^2_{H^{\frac12}(\partial\Omega)}+\int_{S_r^+(x_0)}|\nabla_T u|^2
+\int_{\Omega}|\nabla u|^2+\int_{B_r^+(x_0)}|\tau||\nabla u|\big).
\end{equation}
This yields \eqref{radial_tangent2} and completes the proof. \qed

\section{No energy concentration and oscillation accumulation in boundary neck regions}
\setcounter{equation}{0}
\setcounter{theorem}{0}

This section is devoted to the proof that there is neither energy concentration nor oscillation 
accumulation of the sequence in any boundary neck region, which is defined to be the region either between two consecutive bubbles  or between a boundary bubble and the body region at a boundary point.
 
The crucial step is to show the tangential energy over dyadic boundary annual regions enjoys power decays
with respect to their radius, see \cite{Li-Wang1, Li-Wang2} for related interior estimates. 
To better present this, we need to introduce some notations.

Let $r_0=r_0(\partial\Omega)>0$ be the smallest constant among Lemma 2.2, Theorem 2.3, and Lemma 2.4.
For $0<r<r_0$, $x_0\in\partial\Omega$,  and $t_0>0$, we set a family of dyadic boundary annuals by
$$Q^+(t, t_0; x_0, r)=B^+_{e^{-(t_0-t)} r}(x_0)\setminus B^+_{e^{-(t_0+t)} r}(x_0), \ t>0,$$
and define 
$$E_n(t, t_0; x_0, r)=\int_{Q^+(t, t_0; x_0, r)}|\nabla u_n|^2, \ t>0.$$
Then, by direct calculations, we have that for a.e. $t>0$,
\begin{equation}\label{dE_n}
\frac{d}{dt}E_n(t, t_0; x_0, r)=e^{-(t_0-t)} r\int_{S_{e^{-(t_0-t)}r}^+(x_0)}|\nabla u_n|^2
+e^{-(t_0+t)} r\int_{S_{e^{-(t_0+t)}r}^+(x_0)}|\nabla u_n|^2.
\end{equation}
Now we have
\begin{lemma}\label{neck_energy_decay} For any $\epsilon>0$, 
if $\{u_n\}\subset H^2(\Omega, N)$ is a sequence of approximate harmonic maps
given by Theorem \ref{bubbling1} such that for sufficiently small $\delta>0$ and sufficiently large $R>1$, 
\begin{equation}\label{small_energy}
\sup_{n\ge 1}\sup_{Rr_n\le \tau\le 2\delta}\int_{B_{2\tau}^+(x_n)\setminus B_{\tau}^+(x_n)}|\nabla u_n|^2\le\epsilon^2,
\end{equation}
holds for some $x_n\in\partial\Omega$ and $r_n\rightarrow 0$, 
then there exists $C>0$, independent of $n$, such that 
\begin{equation}\label{energy_decay}
\int_{B_{\frac{\delta}2}^+(x_n)\setminus B_{4Rr_n}^+(x_n)}|\nabla u_n|^2
\le {C}\big(\epsilon+\delta\big),
\end{equation}
and
\begin{equation}\label{osc_decay}
{\rm{osc}}_{B_{\frac{\delta}4}^+(x_n)\setminus B_{8Rr_n}^+(x_n)}u_n
\le C\big(\sqrt{\epsilon}+\sqrt{\delta}\big).
\end{equation}
\end{lemma}
\pf Define a family of radial functions
$$\phi_n(\tau)=\frac{1}{|S_\tau^+(x_n)|}\int_{S_\tau^+(x_n)} u_n\,dH^1, R{r_n}\le\tau\le 2\delta.$$
By the assumption \eqref{small_energy}, we can apply both Lemma \ref{interior} and Theorem \ref{bdry_max1} to conclude that
\begin{eqnarray}\label{oscillation-estimate30}
&&\sup_{x\in B_\delta^+(x_n)\setminus B_{2Rr_n}^+(x_n)}\big|u_n(x)-\phi_n(|x-x_n|)\big|
\le \sup_{2Rr_n\le \tau\le \delta} {\rm{osc}}_{S_\tau^+(x_n)} u_n\nonumber\\
&&\le C\big(\delta\|\tau_n\|_{L^2(B_{2\delta}^+(x_n))}+\sup_{Rr_n\le\tau\le 2\delta}\|\nabla u_n\|_{L^2(B_{2\tau}^+(x_n)\setminus
B_\tau^+(x_n))}+\delta \|u_n-g_n\|_{L^2(T_{2\delta}(x_n))}\big)\nonumber\\
&&\le C(\epsilon+\delta).
\end{eqnarray}
Now choose $r_1=r_1(n)\in (\frac{\delta}2, \delta)$ and $R_2=R_2(n)\in (2Rr_n, 4Rr_n)$
such that
\beq
\label{good_radius}
\begin{cases} 
r_1\int_{S_{r_1}^+(x_n)}|\nabla u_n|^2\le 8\int_{B_\delta^+(x_n)\setminus B_{\frac{\delta}2}^+(x_n)}|\nabla u_n|^2\le
8\epsilon^2,\\
R_2\int_{S_{R_2}^+(x_n)}|\nabla u_n|^2\le 8\int_{B_{4Rr_n}^+(x_n)\setminus B_{2Rr_n}^+(x_n)}|\nabla u_n|^2\le
8\epsilon^2.
\end{cases}
\eeq
Multiplying \eqref{HM_WA00}$_1$ by $u_n(\cdot)-\phi_n(|\cdot-x_n|)$ and integrating over $B_{r_1}^+(x_n)\setminus B_{R_2}^+(x_n)$,
we obtain
\begin{eqnarray*}
&&\int_{B_{r_1}^+(x_n)\setminus B_{R_2}^+(x_n)}\nabla u_n\cdot\nabla(u_n-\phi_n(|\cdot-x_n|)\\
&&=\big(\int_{S_{r_1}^+(x_n)}+\int_{S_{R_2}^+(x_n)}\big)\frac{\partial u_n}{\partial\nu}(u_n-\phi_n(|\cdot-x_n|)\\
&&\ -w\int_{T_{r_1}(x_n)\setminus T_{R_2}(x_n)}(u_n-g_n)(u_n-\phi_n(|\cdot-x_n|)\\
&&+\int_{B_{r_1}^+(x_n)\setminus B_{R_2}^+(x_n)}(A(u_n)(\nabla u_n,\nabla u_n)-\tau_n)(u_n-\phi_n(|\cdot-x_n|)
\end{eqnarray*}
Applying Poincar\'e's inequality and \eqref{good_radius}, we can bound
\begin{eqnarray*}
&&\big|\big(\int_{S_{r_1}^+(x_n)}+\int_{S_{R_2}^+(x_n)}\big)\frac{\partial u_n}{\partial\nu}(u_n-\phi_n(|\cdot-x_n|)\big|\\
&&\le C\big(r_1\int_{S_{r_1}^+(x_n)}|\nabla u_n|^2+R_2\int_{S_{R_1}^+(x_n)}|\nabla u_n|^2\big)\le C\epsilon^2.
\end{eqnarray*}
Using \eqref{oscillation-estimate30}, we can estimate
\begin{eqnarray*}
&&\big|\int_{T_{r_1}(x_n)\setminus T_{R_2}(x_n)}(u_n-g_n)(u_n-\phi_n(|\cdot-x_n|)\big|\\
&&\le C\|u_n-\phi_n(|\cdot-x_n|)\|_{L^\infty(B_{r_1}^+(x_n)\setminus B_{R_2}^+(x_n))}\int_{T_{r_1}(x_n)}|u_n-g_n|
\le C(\epsilon+\delta),
\end{eqnarray*}
and
\begin{eqnarray*}
&&\big|\int_{B_{r_1}^+(x_n)\setminus B_{R_2}^+(x_n)}(A(u_n)(\nabla u_n,\nabla u_n)-\tau_n)(u_n-\phi_n(|\cdot-x_n|))\big|\\
&&\le C\|u_n-\phi_n(|\cdot-x_n|)\|_{L^\infty(B_{r_1}^+(x_n)\setminus B_{R_2}^+(x_n))}\int_{\Omega}(|\nabla u_n|^2+|\tau_n|)
\le C(\epsilon+\delta).
\end{eqnarray*}
While by using polar coordinates we have
\begin{eqnarray}\label{tangential_energy_decay0}
&&\int_{B_{r_1}^+(x_n)\setminus B_{R_2}^+(x_n)}\nabla u_n\cdot\nabla(u_n-\phi_n(|\cdot-x_n|))\nonumber\\
&&=\int_{B_{r_1}^+(x_n)\setminus B_{R_2}^+(x_n)}(|\nabla u_n|^2-|\frac{\partial u_n}{\partial r}|^2)
+ \int_{R_2}^{r_1}\int_{S_r^+(x_n)}
\langle\frac{\partial u_n}{\partial r}, \frac{\partial (u_n-\phi_n(|\cdot-x_n|))}{\partial r}\rangle \, dH^1(\theta)dr
\nonumber\\
&&\ge \int_{B_{r_1}^+(x_n)\setminus B_{R_2}^+(x_n)}(|\nabla u_n|^2-|\frac{\partial u_n}{\partial r}|^2),
\end{eqnarray}
since we have, by H\"older's inequality, that
\begin{eqnarray*}
&&\int_{R_2}^{r_1}\int_{S_r^+(x_n)}
\langle\frac{\partial u_n}{\partial r}, \frac{\partial (u_n-\phi_n(|\cdot-x_n|))}{\partial r}\rangle  dH^1(\theta)dr\\
&&=\int_{R_2}^{r_1}\Big(\int_{S_r^+(x_n)}
|\frac{\partial u_n}{\partial r}|^2-\frac{1}{|S_r^+(x_n)|}\big|\int_{S_r^+(x_n)}
\frac{\partial u_n}{\partial r}\big|^2\Big) dH^1(\theta)dr\ge 0.
\end{eqnarray*}
Thus we obtain
\beq\label{tangent_energy_decay}
\int_{B_{r_1}^+(x_n)\setminus B_{R_2}^+(x_n)}(|\nabla u_n|^2-|\frac{\partial u_n}{\partial r}|^2)
\le C(\epsilon+\delta).
\eeq
On the other hand, by integrating \eqref{radial_tangent2} over $r\in [R_2, r_1]$ and applying
\eqref {tangent_energy_decay}, we have that
\begin{eqnarray}\label{radial_energy_est1}
\int_{B_{r_1}^+(x_n)\setminus B_{R_2}^+(x_n)}\big|\frac{\partial u_n}{\partial r}\big|^2
&\le& C\int_{B_{r_1}^+(x_n)\setminus B_{R_2}^+(x_n)}\big(|\nabla u_n|^2-|\frac{\partial u_n}{\partial r}|^2\big)\nonumber\\
&+& C\Big(1+\|g_n\|^2_{H^\frac12(\partial\Omega)}+\int_{\Omega}|\nabla u_n|^2+
\int_{B_{r_1}^+(x_n)}|\tau_n||\nabla u_n|)\Big) r_1\nonumber\\
&\le&C(\epsilon+\delta).
\end{eqnarray}
Adding \eqref{tangent_energy_decay} with \eqref{radial_energy_est1} yields 
$$\int_{B_{r_1}^+(x_n)\setminus B_{R_2}^+(x_n)}|\nabla u_n|^2\le C(\epsilon+\delta),$$
which gives \eqref{energy_decay}.

To prove \eqref{osc_decay}, we need to perform the above argument in dyadic boundary annuals. 
Let $L(n,\delta)$ be the positive integer $m$ such that
$$4Rr_n e^m\le \frac{\delta}2\le 4Rr_n e^{m+1}, 
\ {\rm{or\ equivalently}}\ L(n,\delta)=\big[\ln (\frac{\delta}{8Rr_n})\big],$$
where $[t]$ denotes the largest integer part of $t$. 

For $1\le t_0\le L(n, \delta)$ and $0\le t\le \min\{t_0, L(n,\delta)-t_0\}$,
multiplying \eqref{HM_WA00}$_1$ by $u_n-\phi_n(|\cdot-x_n|)$ and integrating the resulting equation over $Q^+(t, t_0; x_n, \delta/2)$, we obtain
\begin{eqnarray}\label{energy-estimate30}
&&\int_{Q^+(t, t_0; x_n, \delta/2)}\langle\nabla u_n, \nabla(u_n-\phi_n(|\cdot-x_n|))\rangle\nonumber\\
&&=\int_{\partial Q^+(t, t_0; x_n, \delta/2)}\langle\frac{\partial u_n}{\partial\nu}, u_n-\phi_n(|\cdot-x_n|)\rangle\nonumber\\
&&\ +\int_{Q^+(t, t_0; x_n, \delta/2)}\langle A(u_n)(\nabla u_n,\nabla u_n)-\tau_n, u_n-\phi_n(|\cdot-x_n|)\rangle\nonumber\\
&&=A_n+B_n.
\end{eqnarray}
Observe that, similar to the estimate \eqref{tangential_energy_decay0},  the left hand side can be bounded by
\begin{eqnarray}\label{LHS-estimate}
\int_{Q^+(t, t_0; x_n, \delta/2)}\langle\nabla u_n, \nabla(u_n-\phi_n(|\cdot-x_n|))\rangle
\ge \int_{Q^+(t, t_0; x_n, \delta/2)}\big(|\nabla u_n|^2-|\frac{\partial u_n}{\partial r}|^2\big).
\end{eqnarray}
Write $\partial Q^+(t, t_0; x_n, \delta/2)=\partial_+Q^+(t, t_0; x_n, \delta/2)
\cup\partial_0 Q^+(t, t_0; x_n, \delta/2)$, where
$$\begin{cases}\partial_+Q^+(t, t_0; x_n, \delta/2)=S_{e^{-(t_0-t)}\delta/2}^+(x_n)
\cup S_{e^{-(t_0+t)}\delta/2}^+(x_n),\\
\partial_0 Q^+(t, t_0; x_n, \delta/2)=T_{e^{-(t_0-t)}\delta/2}(x_n)\setminus T_{e^{-(t_0+t)}\delta/2}(x_n).
\end{cases}
$$
Then we can estimate $A_n$ by
\begin{eqnarray*}
A_n&=&\int_{\partial_+ Q^+(t, t_0, x_n, \delta/2)}\langle\frac{\partial u_n}{\partial\nu}, u_n-\phi_n(|\cdot-x_n|)\rangle\\
&&+\int_{\partial_0 Q^+(t, t_0; x_n, \delta/2)}\langle\frac{\partial u_n}{\partial\nu}, u_n-\phi_n(|\cdot-x_n|)\rangle\\
&=& C_n+D_n.
\end{eqnarray*}
Applying \eqref{HM_WA00}$_2$, we can estimate $D_n$ by
\begin{eqnarray}\label{Dn-estimate}
|D_n|&\le& |w_n|\int_{\partial_0 Q^+(t, t_0; x_n, \delta/2)}|u_n-g_n||u_n-\phi_n(|\cdot-x_n|)|\nonumber\\
&\le& C\|u_n-\phi_n(|\cdot-x_n|)\|_{L^\infty(B_{\delta/2}^+(x_n)\setminus B_{2Rr_n}^+(x_n))}
\int_{T_{e^{-(t_0-t)}\delta/2}(x_n)}|u_n-g_n|\nonumber\\
&\le& C(\epsilon+\delta) e^{-(t_0-t)} \delta.
\end{eqnarray}
We can apply the Poincar\'e inequality to estimate $C_n$ as follows:
\begin{eqnarray}\label{Cn-estimate}
&&|C_n|\nonumber\\
&\le& \int_{S_{e^{-(t_0-t)}\delta/2}^+(x_n)}|\frac{\partial u_n}{\partial\nu}||u_n-\phi_n(|\cdot-x_n|)|\nonumber\\
&&+\int_{S_{e^{-(t_0+t)}\delta/2}^+(x_n)}|\frac{\partial u_n}{\partial\nu}||u_n-\phi_n(|\cdot-x_n|)|\nonumber\\
&\le& \big(\int_{S_{e^{-(t_0-t)}\delta/2}^+(x_n)}|\frac{\partial u_n}{\partial\nu}|^2\big)^\frac12 
\big(\int_{S_{e^{-(t_0-t)}\delta/2}^+(x_n)}|u_n-\phi_n(|\cdot-x_n|)|^2\big)^\frac12\nonumber\\
&&\ + \big(\int_{S_{e^{-(t_0+t)}\delta/2}^+(x_n)}|\frac{\partial u_n}{\partial\nu}|^2\big)^\frac12 
\big(\int_{S_{e^{-(t_0+t)}\delta/2}^+(x_n)}|u_n-\phi_n(|\cdot-x_n|)|^2\big)^\frac12\nonumber\\
&\le& Ce^{-(t_0-t)}\delta \big(\int_{S_{e^{-(t_0-t)}\delta/2}^+(x_n)}|\frac{\partial u_n}{\partial\nu}|^2\big)^\frac12 
\big(\int_{S_{e^{-(t_0-t)}\delta/2}^+(x_n)}|\nabla_T u_n|^2\big)^\frac12\nonumber\\
&&\ + Ce^{-(t_0+t)}\delta \big(\int_{S_{e^{-(t_0+t)}\delta/2}^+(x_n)}|\frac{\partial u_n}{\partial\nu}|^2\big)^\frac12 
\big(\int_{S_{e^{-(t_0+t)}\delta/2}^+(x_n)}|\nabla_T u_n|^2\big)^\frac12\nonumber\\
&\le& C\Big(e^{-(t_0-t)}\delta/2 \int_{S_{e^{-(t_0-t)}\delta/2}^+(x_n)}|\nabla u_n|^2
+e^{-(t_0+t)}\delta/2 \int_{S_{e^{-(t_0+t)}\delta/2}^+(x_n)}|\nabla u_n|^2\Big)\nonumber\\
&=& C\frac{d}{dt}E_n(t, t_0; x_n, \delta/2).
\end{eqnarray}
It follows from \eqref{Cn-estimate} and \eqref{Dn-estimate} that
\begin{equation}\label{An-estimate}
|A_n|\le C\Big(\frac{d}{dt}E_n(t, t_0; x_n, \delta/2)+(\epsilon+\delta)e^{-(t_0-t)} \delta\Big).
\end{equation}
Applying \eqref{oscillation-estimate30}, we can estimate $B_n$ by
\begin{eqnarray}\label{Bn-estimate}
|B_n|&\le& C\int_{Q^+(t, t_0; x_n, \delta/2)}(|\nabla u_n|^2+|\tau_n|)
\|u_n-\phi_n(|\cdot-x_n|)\|_{L^\infty(B_{\delta}^+(x_n)\setminus B_{2Rr_n}^+(x_n))}\nonumber\\
&\le& C(\epsilon+\delta) \big(E_n(t, t_0; x_n, \delta/2)
+\|\tau_n\|_{L^2(B_{2\delta}^+(x_n))} e^{-(t_0-t)} \delta\big)\nonumber\\
&\le& C(\epsilon+\delta) \big(E_n(t, t_0; x_n, \delta/2)+ e^{-(t_0-t)} \delta\big).
\end{eqnarray}
Substituting \eqref{An-estimate} and \eqref{Bn-estimate} into \eqref{energy-estimate30} and applying
\eqref{LHS-estimate}, we arrive at
\begin{eqnarray}\label{energy-estimate31}
\int_{Q^+(t, t_0; x_n, \delta/2)}\big(|\nabla u_n|^2-|\frac{\partial u_n}{\partial r}|^2\big)
&\le& C(\epsilon+\delta)E_n(t, t_0; x_n, {\delta}/2)\nonumber\\
&+& C\frac{d}{dt}E_n(t, t_0; x_n, {\delta}/2)
+C(\epsilon+\delta) e^{-(t_0-t)} \delta.
\end{eqnarray}
Now we apply \eqref{radial_tangent2} and \eqref{energy-estimate31} to get
\begin{eqnarray}\label{energy-estimate32}
\int_{Q^+(t, t_0; x_n, {\delta}/2)}\big|\frac{\partial u_n}{\partial r}\big|^2
&\le& C\int_{Q^+(t, t_0; x_n, {\delta}/2)}\big(|\nabla u_n|^2-|\frac{\partial u_n}{\partial r}|^2\big)\nonumber\\
&+& C\Big(1+\|g_n\|^2_{H^\frac12(\partial\Omega)}+\int_{\Omega}|\nabla u_n|^2+
\int_{B_{2\delta}^+(x_n)}|\tau_n||\nabla u_n|)\Big) e^{-(t_0-t)}\delta
\nonumber\\
&\le&C\int_{Q^+(t, t_0; x_n, {\delta}/2)}\big(|\nabla u_n|^2-|\frac{\partial u_n}{\partial r}|^2\big) +C(1+\sqrt{\epsilon})e^{-(t_0-t)}\delta\nonumber\\
&\le& C\big((\epsilon+\delta)E_n(t, t_0; x_n, {\delta}/2)+\frac{d}{dt}E_n(t, t_0; x_n, {\delta}/2)
+e^{-(t_0-t)} {\delta}\big).
\end{eqnarray}
Adding \eqref{energy-estimate31} with \eqref{energy-estimate32}, we arrive at
\begin{equation}\label{energy-estimate33}
E_n(t, t_0; x_n, {\delta}/2)\le C(\epsilon+\delta)E_n(t, t_0; x_n, {\delta}/2)
+C\frac{d}{dt}E_n(t, t_0; x_n, {\delta}/2)+Ce^{-(t_0-t)} \delta.
\end{equation}
By choosing sufficiently small $\epsilon$ and $\delta$ so that 
$C(\epsilon+\delta)\le \frac12$, this implies that there exists $0<c<1$ such that 
\begin{equation}\label{energy-estimate34}
cE_n(t, t_0; x_n, {\delta}/2)\le \frac{d}{dt}E_n(t, t_0; x_n, {\delta}/2)+e^{-(t_0-t)} \delta, 
\end{equation}
and hence
\begin{equation}\label{ODE}
\frac{d}{dt}\big(e^{-ct}E_n(t, t_0; x_n, {\delta}/2)\big)\ge -e^{-t_0+(1-c)t}{\delta}.
\end{equation}
Set $1\le t_0=i\le L(n,\delta)$ and integrate \eqref{ODE} over $t\in [1, m(i, n)]$, where
$m(i, n)=\min\big\{i, L(n,\delta)-i\big\}$, we obtain that
\begin{eqnarray}\label{decay1}
\int_{B^+_{e^{-i+1}\delta/2}(x_n)
\setminus B^+_{e^{-i-1}\delta/2}(x_n)}|\nabla u_n|^2
&=&E_n(1, i; x_n, \delta/2)\nonumber\\
&\le& C(\int_{B_{\frac{\delta}2}^+(x_n)\setminus B_{4Rr_n}^+(x_n)}|\nabla u_n|^2+\delta)e^{-cm(i, n)}
\nonumber\\
&\le&C(\epsilon+\delta)e^{-cm(i,n)},
\end{eqnarray}
where we have used \eqref{energy_decay} in the last step.
It follows from \eqref{decay1} that
\begin{eqnarray*}
\int_{B_{\frac{\delta}4}^+(x_n)\setminus B_{8Rr_n}^+(x_n)}\frac{|\nabla u_n|}{|x-x_n|}
&\le&\sum_{i=1}^{m(1,n)}
\int_{B_{e^{-i+1}\delta/2}^+(x_n)\setminus B_{e^{-i-1}\delta/2}^+(x_n)}\frac{|\nabla u_n|}{|x-x_n|}\\
&\le& C \sum_{i=1}^{m(1,n)}\big(\int_{B^+_{e^{-i+1}\delta/2}(x_n)
\setminus B^+_{e^{-i-1}\delta/2}(x_n)}|\nabla u_n|^2\big)^\frac12\\
&\le &C\sum_{i=1}^{m(1,n)}\sqrt{E_n(1, i; x_n, \delta/2)}\\
&\le& C(\sqrt{\epsilon}+\sqrt{\delta})\sum_{i=1}^{m(1,n)}e^{-\frac{cm(i,n)}2}
\le C(\sqrt{\epsilon}+\sqrt{\delta}).
\end{eqnarray*}
On the other hand, it follows from direct calculations and \eqref {oscillation-estimate30} that 
\begin{eqnarray*}
{\rm{osc}}_{B_{\frac{\delta}4}^+(x_n)\setminus B_{8Rr_n}^+(x_n)}u_n
&\le& C\Big(\int_{B_{\frac{\delta}4}^+(x_n)\setminus B_{8Rr_n}^+(x_n)}\frac{|\nabla u_n|}{|x-x_n|}
+\sup_{2Rr_n\le r\le \delta} {\rm{osc}}_{S_r^+(x_n)} u_n\Big)\\
&\le& C\big(\sqrt{\epsilon}+\sqrt{\delta}\big).
\end{eqnarray*}
This implies \eqref{osc_decay}. The proof is now complete. \qed


\section{Removable isolated singularity at the boundary}
\setcounter{equation}{0}
\setcounter{theorem}{0}
In order to show the weak limit $u$ in Theorem \ref{bubbling1} belongs to $H^2(\Omega, N)$, we need to establish
the removability of an isolated singularity, both in the interior and on the boundary of $\Omega$,  for an approximate harmonic map under weak anchoring condition, with tension field $\tau\in L^2(\Omega, T_uN)$.  The following removability of an interior isolated singularity has been known before. 
\begin{lemma}\label{remove_interior} Assume $u\in H^1(\Omega, N)$ is an approximate harmonic map
with tension field $\tau\in L^2(\Omega, T_uN)$. If $u\in H^2_{\rm{loc}}(B_{r_0}(x_0)\setminus\{x_0\}, N)$ 
for some $B_{r_0}(x_0)\subset\Omega$, then $u\in H^2(B_{r_0}(x_0), N)$.
\end{lemma}
\pf See \cite{Ding-Tian}. \qed

\medskip
Now we want to prove the following result on the removability of a boundary isolated singularity.

\begin{lemma} \label{remove_bdry}
For $x_0\in\partial\Omega$ and $r_0>0$, assume 
that $u\in H^1(B_{r_0}^+(x_0), N)$ solves 
\begin{equation}\label{HM_WA40}
\begin{cases}
\Delta u+A(u)(\nabla u,\nabla u) =\tau & \ {\rm{in}}\ B_{r_0}^+(x_0),\\
\frac{\partial u}{\partial \nu}+w \mathbb P(u)(u-g)=0 & \ {\rm{on}}\  T_{r_0}(x_0),
\end{cases}
\end{equation}
for some $\tau\in L^2(B_{r_0}^+(x_0), T_uN)$, $w>0$,
and $g\in H^{\frac12}(T_{r_0}(x_0), N)$. 
If $u\in H^2_{\rm{loc}}(B_{r_0}^+(x_0)\setminus\{x_0\}, N)$,
then $u\in H^2(B_{r_0}^+(x_0), N)$.
\end{lemma}
\pf For simplicity, we assume $x_0=0$, $r_0=1$, and $\Omega=\mathbb R^2_+:=\big\{x=(x_1,x_2)\in\R^2:
x_2>0\big\}$. Applying the same argument as Lemma \ref{neck_energy_decay}, we can prove that there exists
a sufficiently small $r_1>0$ and $\alpha\in (0,1)$ such that
\beq\label{energy_power_decay}
\int_{B_r^+(0)}|\nabla u|^2\le Cr^{\alpha}, \ \forall\  0<r\le r_1.
\eeq
Here we sketch the proof of \eqref{energy_power_decay}. Since $u\in H^1(B_1^+(0))$, there exists 
$0<r_1\le \frac12$ such that 
$$\int_{B_{r_1}^+(x)}|\nabla u|^2\le \epsilon_0^2, \ \forall\ x\in B_{\frac12}^+(0).$$
For any $x\in B_{r_1}^+(0)\setminus\{0\}$, we then have
\beq\label{dadic_small}
\int_{B_{\frac{|x|}2}(x)\cap\R^2_+}|\nabla u|^2\le \epsilon_0^2,
\eeq
and, since $u\in H^2_{\rm{loc}}(B_1^+(0)\setminus\{0\})$ and $B_{\frac{|x|}2}(x)\cap \R^2_+\Subset B_1^+(0)\setminus\{0\}$, we have $u\in H^2(B_{\frac{|x|}2}(x)\cap\R^2_+)$.
Thus we can apply Lemma \ref{interior} and Theorem \ref{bdry_max1} to conclude that
$u\in C(B_{\frac{|x|}4}(x)\cap\R^2_+)$ and
$$\big(\underset{B_{\frac{|x|}4}(x)\cap\R^2_+}{\rm{osc}} u\big)^2
\le C\big(\int_{B_{\frac{|x|}2}(x)\cap\R^2_+}(|\nabla u|^2+|x|^2|\tau|^2)+w^2|x|\int_{B_{\frac{|x|}2}(x)\cap\partial\R^2_+}
|u-g|^2\big).
$$
Therefore we obtain that
\beq\label{osc_estimate3}
\sup_{0<r\le r_1}\underset{S_{r}^+(0)}{\rm{osc}} u
\le C(\epsilon_0+\|\tau\|_{L^2(B_1^+(0))} r_1+\|u-g\|_{L^2(T_1(0))} r_1^\frac12)\le C(\epsilon_0+r_1^\frac12).
\eeq
Set
$$\phi(r)=\frac{1}{|S_r^+(0)|}\int_{S_r^+(0)} u\,dH^1, \ 0 <r\le r_1.$$
Then we have
\beq\label{osc_estimate4}
\max_{x\in B_{r_1}^+(0)\setminus\{0\}} |u(x)-\phi(|x|)|\le \sup_{0<r\le r_1}\underset{S_{r}^+(0)}{\rm{osc}} u
\le C(\epsilon_0+r_1^\frac12).
\eeq
For $0<s<r\le r_1$, multiplying \eqref{HM_WA40}$_1$ by $u-\phi$ and integrating over $B_r^+(0)\setminus B_s^+(0)$,
we obtain that
\begin{eqnarray*}
&&\int_{B_r^+(0)\setminus B_s^+(0)}\nabla u\cdot\nabla(u-\phi)\\
&&=\int_{\partial(B_r^+(0)\setminus B_s^+(0))}\frac{\partial u}{\partial\nu}\cdot (u-\phi)
+\int_{B_r^+(0)\setminus B_s^+(0)}(A(u)(\nabla u, \nabla u)-\tau)(u-\phi)\\
&&=\big(\int_{S_r^+(0)}+\int_{S_s^+(0)}\big)\frac{\partial u}{\partial\nu}\cdot (u-\phi)-
w\int_{T_r(0)\setminus T_s(0)}(u-g)(u-\phi)\\
&&\ \ \ +\int_{B_r^+(0)\setminus B_s^+(0)}(A(u)(\nabla u, \nabla u)-\tau)(u-\phi).
\end{eqnarray*}
Choosing a sequence $s\rightarrow 0$  so that
$$s\int_{S_s^+(0)}|\nabla u|^2\rightarrow 0$$
and hence, by sending $s\rightarrow 0$, we arrive at
\begin{eqnarray*}
&&\int_{B_r^+(0)}\nabla u\cdot\nabla(u-\phi)\\
&&=\int_{S_r^+(0)}\frac{\partial u}{\partial\nu}\cdot (u-\phi)-w\int_{T_r(0)}(u-g)(u-\phi)\\
&&+\int_{B_r^+(0)}(A(u)(\nabla u, \nabla u)-\tau)(u-\phi).
\end{eqnarray*}
As in the proof of Lemma \ref{neck_energy_decay}, we can bound
$$\big|\int_{S_r^+(0)}\frac{\partial u}{\partial\nu}\cdot (u-\phi)\big|\le Cr\int_{S_r^+(0)}|\nabla u|^2,$$
$$\big|\int_{T_r(0)}(u-g)(u-\phi)\big|\le C\|u-\phi\|_{L^\infty(S_r^+(0))}\int_{T_r(0)}|u-g|\le Cr,$$
$$\big|\int_{B_r^+(0)}(A(u)(\nabla u, \nabla u)-\tau)(u-\phi)\big|
\le C(\epsilon_0+r^\frac12)\big(\int_{B_r^+(0)}|\nabla u|^2+r\|\tau\|_{L^2(B_r^+(0))}\big),
$$
and
$$
\int_{B_r^+(0)}\nabla u\cdot\nabla(u-\phi)\ge \int_{B_r^+(0)}\big(|\nabla u|^2-|\frac{\partial u}{\partial r}|^2\big).
$$
Hence we obtain 
\beq\label{tangential_energy_est}
\int_{B_r^+(0)}\big(|\nabla u|^2-|\frac{\partial u}{\partial r}|^2\big)
\le Cr\int_{S_r^+(0)}|\nabla u|^2+C(\epsilon_0+r^\frac12)\big(\int_{B_r^+(0)}|\nabla u|^2+r\|\tau\|_{L^2(B_r^+(0))}\big)
+Cr.
\eeq

Now we claim that \eqref{radial_tangent2} of Lemma \ref{radial_tangent1} also holds 
for an approximate harmonic map $u\in H^2_{\rm{loc}}(B_1^+(0)\setminus\{0\})$. In fact, this can be
achieved by applying the same argument as in \eqref{pohozaev} of Lemma \ref{radial_tangent1} with the integration
domain $B_r^+(0)$ replaced by $B_r^+(0)\setminus B_{s_k}^+(0)$, where $s_k\rightarrow 0$ is chosen
so that
$$s_k\int_{S_{s_k}^+(0)}|\nabla u|^2\rightarrow 0.$$ 
Hence we also have
\begin{eqnarray}\label{radial_tangent3}
\int_{B_r^+(0)}\big|\frac{\partial u}{\partial r}\big|^2
&\le& C\int_{B_r^+(0)}\big(|\nabla u|^2-|\frac{\partial u}{\partial r}|^2\big)\nonumber\\
&+& C\Big(1+\|g\|^2_{H^\frac12(T_1(0))}+\int_{B_1^+(0)}|\nabla u|^2+
\int_{B_{1}^+(0)}|\tau||\nabla u|)\Big) r.
\end{eqnarray}
Adding \eqref{tangential_energy_est} with \eqref{radial_tangent3}, and choosing sufficiently small
$\epsilon_0>0$ and $r_1>0$,  we obtain that
\beq\label{morrey_decay}
\int_{B_r^+(0)}|\nabla u|^2\le Cr\int_{S_r^+(0)}|\nabla u|^2+Cr,\ \forall\ 0<r\le r_1.
\eeq
This, after integrating over $r$, yields that for some $0<\alpha<1$, it holds
$$\int_{B_r^+(0)}|\nabla u|^2\le r^\alpha
\Big(r_1^{-\alpha}\int_{B_{r_1}^+(0)}|\nabla u|^2+\frac{C}{1-\alpha}r_1^{1-\alpha}\Big), \ \forall\  0<r\le r_1.
$$
Hence \eqref{energy_power_decay} holds.

Next we claim that there exist $p>1$ and $C=C(p)>0$ such that 
\beq\label{Lp-bound}
\int_{B_{\frac{r_1}2}^+(0)}|\nabla u|^{2p}\le C.
\eeq
It follows from \eqref{dadic_small} that we can apply Lemma \ref{interior},
Theorem \ref{bdry_max1}, and  \eqref{energy_power_decay} to show that
for any $0<r\le \frac{r_1}4$,
\begin{eqnarray*}
\big(r^2\int_{B_{2r}^+(0)\setminus B_r^+(0)}|\nabla u|^4\big)^\frac14
&\le& C\big(\|\nabla u\|_{L^2(B_{2r}^+(0)\setminus B_r^+(0))}+r\|\nabla^2 u\|_{L^2(B_{2r}^+(0)\setminus B_r^+(0))}\big)\\
&\le& C\big(\|\nabla u\|_{L^2(B_{4r}^+(0)\setminus B_{\frac{r}2}^+(0))}+r\|\tau\|_{L^2(B_{4r}^+(0))}
+r\|g\|_{H^\frac12(T_{4r}(0))}\big)\\
&\le& Cr^{\frac{\alpha}2}.
\end{eqnarray*}
By H\"older's inequality, this implies that for any $1<p\le 2$,
$$\int_{B_{2r}^+(0)\setminus B_r^+(0)}|\nabla u|^{2p}\le Cr^{\alpha p+2-2p}, \ 0< r\le \frac{r_1}4.$$
Hence, after choosing $1<p<\frac{2}{2-\alpha}$, we have that
\begin{eqnarray*}
\int_{B_{\frac{r_1}2}^+(0)}|\nabla u|^{2p}
&=&\sum_{i=0}^\infty \int_{B_{2^{-i-1}r_1}^+(0)\setminus B_{2^{-i-2} r_1}^+(0)}|\nabla u|^{2p}\\
&\le& Cr_1^{\alpha p+2-2p}\sum_{i=0}^\infty 2^{-(i+2)(\alpha p+2-2p)}\le Cr_1^{\alpha p+2-2p}.
\end{eqnarray*}
This gives \eqref{Lp-bound}. 

Using \eqref{Lp-bound}, we can apply the standard $W^{2,p}$-theory to \eqref{HM_WA40} to conclude
that $u\in W^{2,p}(B_{\frac{r_1}3}^+(0))$. By Sobolev's embedding theorem, we then have
that $\nabla u\in L^{\frac{2p}{2-p}}(B_{\frac{r_1}3}^+(0))$. Hence by applying the standard $W^{2,p}$-theory
again we conclude that $u\in W^{2, \frac{p}{2-p}}(B_{\frac{r_1}4}^+(0))$. This, along with the bootstrapping argument,
can eventually imply that $u\in H^{2}(B_{\frac{r_1}2}^+(0))$.  \qed
\section{Proof of Theorem \ref{bubbling1}}
\setcounter{equation}{0}
\setcounter{theorem}{0}
In this section, based on the analysis from \S2 to \S4, we will provide the proof of Theorem \ref{bubbling1}.

\medskip
\noindent{\bf Proof of Theorem} \ref{bubbling1}. Let $\epsilon_0>0$ be the smaller constant given
by Lemma \ref{interior} and Theorem \ref{bdry_max1}. Define the concentration set
$$\Sigma=\bigcap_{r>0}\big\{x\in\overline\Omega: \liminf_{n\rightarrow\infty}\int_{B_r(x)\cap\Omega}|\nabla u_n|^2\ge\epsilon_0^2\big\}.$$
It is well-known that $\Sigma$ is a finite set of $m_0$ points, with
$m_0\le \frac{E_0}{\epsilon_0^2}$ and $E_0=\sup_{n\ge 1}\int_{\Omega}|\nabla u_n|^2$.
It follows from Lemma \ref{interior} and Theorem \ref{bdry_max1} that 
$$u_n\rightarrow u \ {\rm{in}}\  H^1_{\rm{loc}}(\overline\Omega\setminus\Sigma)\cap C^0_{\rm{loc}}(\overline\Omega\setminus\Sigma).$$
Moreover, by Lemma \ref{remove_interior} and Lemma \ref{remove_bdry}, we know that $u\in H^2(\Omega, N)$
is an approximate harmonic map with tension field $\tau$ and weak anchoring condition $g$ and anchroing
strength parameter $w$. 

Set
$$\Sigma=\big\{x_1,\cdots, x_{m_0}\big\}=\Sigma_1\cup \Sigma_2=
\big\{x_1,\cdots, x_{i_0}\big\}\cup\big\{x_{i_0+1},\cdots, x_{m_0}\big\},$$
for some $1\le i_0\le m_0$, where $\Sigma_1=\big\{x_1,\cdots, x_{i_0}\big\}=\Sigma\cap\Omega$
and $\Sigma_2=\big\{x_{i_0+1},\cdots, x_{m_0}\big\}=\Sigma\cap\partial\Omega$.
Define
$$\delta_1=\frac12{\rm{dist}}(\Sigma_1, \partial\Omega)>0,
\ r(\Sigma_2)=\frac12\inf\big\{|x-y|: \ x, y\in\Sigma_2, \ x\not=y\big\}>0.$$
Then, from the previous works (see, e.g., \cite{Ding-Tian}) on the interior bubbling of approximate harmonic maps with $L^2$-tension fields, 
we know that after exhausting all possible bubbles generated by the concentration set $\Sigma_1$, we can find
a positive integer $m_1$, blowing up points $\{x_n^i\}_{i=1}^{m_1}\subset \Omega_{\delta_1}=\big\{x\in\Omega: {\rm{dist}}(x,\partial\Omega)>\delta_1\big\}$, blowing up scales $\{r_n^i\}_{i=1}^{m_1}\subset\R_+$ with
$\lim_{n\rightarrow\infty} r_n^i=0$ for $1\le i\le m_1$, and nontrivial harmonic maps $\{\omega_i\}_{i=1}^{m_1}\subset C^\infty(\mathbb S^2, N)$ such that
\beq\label{interior-bubble-scales}
\lim_{n\rightarrow\infty}\big\{\frac{r_n^i}{r_n^j},\ \frac{r_n^j}{r_n^i}, \frac{|x_n^i-x_n^j|}{r_n^i+r_n^j}\big\}=\infty,
\ 1\le i<j\le m_1,
\eeq
\beq\label{interior-energy_id}
\lim_{n\rightarrow\infty}\int_{\Omega_{\delta_1}}|\nabla u_n|^2
=\int_{\Omega_{\delta_1}}|\nabla u|^2+\sum_{i=1}^{m_1}\int_{\mathbb S^2}|\nabla\omega_i|^2,
\eeq
and
\beq\label{interior-oscillation_conv}
\lim_{n\rightarrow\infty}\big\|u_n-u-\sum_{i=1}^{m_1}\big(\omega_i(\frac{\cdot-x_n^i}{r_n^i})-\omega_i(\infty)\big)\big\|_{L^\infty(\Omega_{\delta_1})}=0.
\eeq

Now we need to carry out the blowing up analysis near the boundary concentration set $\Sigma_2$.
For $x_j\in\Sigma_2$, $i_0+1\le j\le m_0$, there exist $0<r_n^j\le r(\Sigma_2)$ and $x_n^j\in B_{r(\Sigma_2)}^+(x_j)$
such that
$$\int_{B_{r_n^j}(x_n^j)\cap\Omega}|\nabla u_n|^2
=\max\big\{\int_{B_{r_n^j}(x)\cap\Omega}|\nabla u_n|^2: x\in B_{r(\Sigma_2)}^+(x_j)\big\}=\frac{\epsilon_0^2}{10^2}.$$

It is readily seen that $x_n^j\rightarrow x_j$ and $r_n^j\rightarrow 0$. Define the blowing up sequence
$$v_n^j(x)=u_n(x_n^j+r_n^j x): \Omega_n^j=(r_n^j)^{-1}(B_{r(\Sigma_2)}^+(x_j)\setminus\{x_j\})\to N.$$
Then $v_n^j$ satisfies\\
i) $$\begin{cases}
\Delta v_n^j+A(v_n^j)(\nabla v_n^j, \nabla v_n^j)=\tau_n^j,\ {\rm{in}}\ \Omega_n^j, \\
\frac{\partial v_n^j}{\partial\nu}+w_n r_n^j \mathbb P(v_n^j)(v_n^j-g_n^j)=0, \ {\rm{on}}\ \partial_0\Omega_n^j,
\end{cases}
$$
where $$\tau(v_n^j)(x)=(r_n^1)^2\tau(u_n)(x_n^1+r_n^1 x), g_n^j(x)=g_n(x_n^j+r_n^j x) 
\ {\rm{for}}\ x\in \Omega_n^j, \, \partial_0\Omega_n^j=(r_n^j)^{-1}(T_{r(\Sigma_2)}(x_j)\setminus\{x_j\}).$$
ii) $\displaystyle\int_{\Omega_n^j}|\nabla v_n^j|^2\le \int_{\Sigma}|\nabla u_n|^2\le E_0.$\\
iii) 
$\displaystyle\int_{B_1(x)\cap\Omega_n^j}|\nabla v_n^j|^2\le \int_{B_1^+(0)}|\nabla v_n^j|^2=\frac{\epsilon_0^2}{10^2}, \ \forall\ x\in\Omega_n^j.$ Hence
 $$\int_{B_{10}(x)\cap\Omega_n^j}|\nabla v_n^j|^2\le \epsilon_0^2, \ \forall \ x\in\Omega_n^j.$$
We now divide the argument into two possible cases:\\
a) $\displaystyle\lim_{n\rightarrow\infty}\frac{{\rm{dist}}(x_n^j, \partial\Omega)}{r_n^j}=\infty$
so that $\Omega_n^j\rightarrow \R^2$ as $n\rightarrow\infty$: Applying Lemma \ref{interior}, we conclude that
there exists a nontrivial harmonic map $\omega_j\in C^\infty\cap H^1(\R^2, N)$ such that
$v_n^j\rightarrow \omega_j$ in $C^0\cap H^1_{\rm{loc}}(\R^2)$. It is well-known that $\omega_j$ can be lifted
into a nontrivial harmonic map from $\mathbb S^2$ to $N$.\\
b) $\displaystyle\lim_{n\rightarrow\infty}\frac{{\rm{dist}}(x_n^j, \partial\Omega)}{r_n^j}=a$ for some $0\le a<\infty$
so that $\Omega_n^j\rightarrow \R^2_{-a}=\big\{x=(x_1,x_2)\in\R^2: x_2\ge -a\big\}$ as $n\rightarrow\infty$:
Applying Lemma \ref{interior} and Theorem \ref{bdry_max1}, we conclude that
there exists a nontrivial harmonic map $\omega_j\in C^\infty\cap H^1(\R^2_{-a}, N)$, with
$\frac{\partial\omega_j}{\partial\nu}=0$ on $\partial\R^2_{-a}$, such that
$v_n^j\rightarrow \omega_j$ in $C^0\cap H^1_{\rm{loc}}(\R^2_{-a})$. However,
it is well-known that any harmonic map $\omega_j\in H^1\cap C^\infty(\R^2_{-a}, N)$, 
with $\frac{\partial\omega_j}{\partial\nu}=0$ on $\partial\R^2_{-a}$, must be constant. We get a desired
contradiction. Thus the case b) doesn't occur. 

Repeating this process for each $x_j\in\Sigma_2$, we can find a positive integer $m_2>m_1$, 
all possible bubbles $\{\omega_j\}$, $m_1+1\le j\le m_2$, generated by $\Sigma_2$, and sequences
of blowing up points and scales $\{x_n^j\}\subset\overline\Omega\setminus\Omega_{\delta_1}$, with
${\rm{dist}}(x_n^j, \Sigma_2)\rightarrow 0$, and $r_n^j\rightarrow 0$, such that
$$u_n(x_n^j+r_n^j x)\rightarrow \omega_j \ {\rm{in}}\ H^1\cap C^0_{\rm{loc}}(\R^2), \ m_1+1\le j\le m_2.$$
Reasoning as in the interior case above, we can see that the property \eqref{bubble-scales} also holds
for $m_1+1\le i<j\le m_2$. This, combined with \eqref{interior-bubble-scales}, yields \eqref{bubble-scales} holds for all $1\le i<j\le m_2$.

From \eqref{interior-energy_id} and \eqref{interior-oscillation_conv},  we see that 
in order to prove \eqref{energy_id} and \eqref{oscillation_conv}, it suffices to show 
\beq
\label{bdry_energy_id}
\lim_{n\rightarrow\infty}\int_{\Omega\setminus\Omega_{\delta_1}}|\nabla u_n|^2
=\int_{\Omega\setminus\Omega_{\delta_1}}|\nabla u|^2+\sum_{i=m_1+1}^{m_2}\int_{\mathbb S^2}|\nabla\omega_i|^2,
\eeq
and
\beq
\label{bdry_oscillation_conv}
\lim_{n\rightarrow\infty}\big\|u_n-u-\sum_{m_1+1}^{m_2}\big(\omega_i(\frac{\cdot-x_n^i}{r_n^i})-\omega_i(\infty)\big)\big\|_{L^\infty(\Omega\setminus\Omega_{\delta_1})}=0.
\eeq
It is well-known that, by the induction scheme of \cite{Ding-Tian}, we can further assume that $m:=m_2=m_1+1$, 
$\Sigma_2=\{x_{m}\}$ consists of one point, and there is only one bubble $\omega_{m}$ generated at $x_{m}$.
As in \cite{Ding-Tian}, this assumption yields that for any $\epsilon>0$, there exist sufficiently large $R>1$
and sufficiently small $0<\delta<\delta_1$ such that
\beq\label{necksmall1}
\sup_{Rr_n^m\le r\le \delta}\int_{(B_{2r}(x_n^m)\setminus B_r(x_n^m))\cap\Omega}|\nabla u_n|^2\le\epsilon^2,
\ \forall\ n\ge 1.
\eeq
Set $d_n^m={\rm{dist}}(x_n^{m}, \partial\Omega)>0$. Then $d_n^m\rightarrow 0$ and 
$\lambda_n=\frac{d_n^m}{r_n^{m}}\rightarrow \infty$.
Let $\{y_n^{m}\}\subset\partial\Omega$ be such that $d_n^m=|x_n^{m}-y_n^{m}|$. Then  we  have
that $y_n^{m}\rightarrow x_{m}$ as $n\rightarrow\infty$.

It follows from the proof of \eqref{interior-energy_id} that it holds
\beq\label{interior_energy_id1}
\lim_{R\rightarrow\infty}\lim_{n\rightarrow\infty}\int_{B_{d_n^m}(x_n^m)\setminus B_{Rr_n^m}(x_n^m)}|\nabla u_n|^2=0,
\eeq
and
\beq\label{interior_oscillation_conv1}
\lim_{R\rightarrow\infty}\lim_{n\rightarrow\infty}\underset{B_{d_n^m}(x_n^m)\setminus B_{Rr_n^m}(x_n^m)}
{\rm{osc}} u_n=0.
\eeq
Observe that
\begin{eqnarray}\label{inclusion1}
\big(B_{\frac{\delta}4}(x_n^m)\cap\Omega\big)\setminus B_{d_n^m}(x_n^m)
= \big((B_{\frac{\delta}4}(x_n^m)\setminus B_{4d_n^m}(x_n^m))\cap\Omega\big)
\cup \big((B_{4d_n^m}(x_n^m)\cap\Omega)\setminus B_{d_n^m}(x_n^m)\big),
\end{eqnarray}
and
\begin{eqnarray}\label{inclusion2}
(B_{\frac{\delta}4}(x_n^m)\setminus B_{4d_n^m}(x_n^m))\cap\Omega
\subset B_{\frac{\delta}2}^+(y_n^m)\setminus B_{2d_n^m}^+(y_n^m)
\subset (B_{\delta}(x_n^m)\cap\Omega)\setminus B_{d_n^m}(x_n^m).
\end{eqnarray}
It follows from \eqref{necksmall1} that
\beq\label{neck_energy5}
\int_{(B_{4d_n^m}(x_n^m)\cap\Omega)\setminus B_{d_n^m}(x_n^m)}|\nabla u_n|^2\le 4\epsilon^2,
\eeq
and
\beq\label{neck_small1}
\sup_{d_n^m\le r\le {\delta}}\int_{B_{2r}^+(y_n^m)\setminus B_r^+(y_n^m)}|\nabla u_n|^2\le \epsilon^2.
\eeq
Thus we can apply Lemma \ref {neck_energy_decay} to obtain that
\beq\label{neck_energy6}
\int_{B^+_{\frac{\delta}2}(y_n^m)\setminus B_{2d_n^m}^+(y_n^m)}|\nabla u_n|^2\le C(\epsilon+\delta),
\eeq
and
\beq\label{neck_osc5}
\underset{B_{\frac{\delta}2}(y_n^m)\setminus B_{2d_n^m}^+(y_n^m)}{\rm{osc}}u_n
\le C(\sqrt\epsilon+\sqrt\delta).
\eeq
It follows from \eqref{neck_energy5} and \eqref{neck_energy6}, and the inclusions \eqref{inclusion1} and
\eqref{inclusion2} that we have
\beq\label{neck_energy7}
\int_{(B_{\frac{\delta}4}(x_n^m)\setminus B_{4d_n^m}(x_n^m))\cap\Omega}|\nabla u_n|^2\le C(\epsilon+\delta),
\eeq
which, combined with \eqref{interior_energy_id1}, implies
\beq\label{neck_energy8}
\lim_{R\rightarrow\infty}\lim_{n\rightarrow\infty}\int_{(B_{\frac{\delta}4}(x_n^m)\cap\Omega)\setminus B_{Rr_n^m}(x_n^m)}|\nabla u_n|^2\le C(\epsilon+\delta).
\eeq
This yields \eqref{bdry_energy_id} after sending $\delta\rightarrow 0$, since $\epsilon>0$ is arbitrarily small.

From the assumption \eqref{neck_small1}, we can apply Lemma \ref{interior} and Theorem \ref{bdry_max1}
to obtain that
\beq\label{neck_osc6}
\underset{B_{2d_n^m}^+(y_n^m)\setminus B_{d_n^m}^+(y_n^m)}{\rm{osc}}u_n\le C(\sqrt\epsilon+\sqrt\delta).
\eeq
Combining \eqref{interior_oscillation_conv1} together with \eqref{neck_osc5} and \eqref{neck_osc6} implies
$$\lim_{R\rightarrow \infty}\lim_{n\rightarrow\infty}\underset{(B_{\frac{\delta}4}(x_n^m)\cap\Omega)\setminus B_{Rr_n^m}(x_n^m)}{\rm{osc}}u_n\le C(\sqrt\epsilon+\sqrt\delta),$$
which yields \eqref{bdry_oscillation_conv} after sending $\delta\rightarrow 0$, since $\epsilon>0$ is arbitrarily small.
The proof of Theorem \ref{bubbling1} is now complete. \qed

\section{Estimates of approximate harmonic maps under Dirichlet conditions}
\setcounter{equation}{0}
\setcounter{theorem}{0}

This section is devoted to the apriori estimates of approximate harmonic maps under Dirichlet conditions. 
First we recall an interior H\"older continuity estimate of approximate harmonic maps $u$ with tension fields
$\tau$ in $\mathcal M^{1,a}(\Omega)$ for some $1< a<2$, which was proved by Wang \cite{wang16}.

\begin{lemma} \label{interior1} There exists $\epsilon_0>0$ such that
if $u\in W^{1,2}(\Omega, N)$ is an approximate harmonic map, with tension field 
$\tau(u)\in M^{1,a}(\Omega)$ for some $1<a<2$, which
satisfies, for some $B_{2r_0}(x_0)\subset\Omega$,
\beq\label{small_cond1}
\displaystyle\int_{B_{2r_0}(x_0)} |\nabla u|^2\le\epsilon_0^2,
\eeq 
then  $u\in C^{2-a}(B_{r_0}(x_0), N)$, and
\begin{equation}\label{holder-estimate}
\big[u\big]_{C^{2-a}(B_{r_0}(x_0))}\le C(a) \big(\epsilon_0+\big\|\tau\big\|_{M^{1,a}(B_{2r_0}(x_0))}\big).
\end{equation}
\end{lemma}
\pf See Wang \cite{wang16} Lemma 2.3. \qed

\medskip
Similar to the proof of Theorem 1.2, in order to deal with the strong anchoring condition, we need to
establish a boundary estimate analogous to Lemma \ref{interior1}.  More precisely, we have
\begin{theorem}\label{boundary1}
There exist  $\ve_0>0$ and $r_0=r_0(\partial\Omega)>0$
such that if $u\in H^{1}(\Omega,N)$ is an approximate harmonic map 
under the Dirichlet condition:
\beq\notag
\left\{
\ba{rll}
\Delta u+A(u)(\nabla u,\nabla u)=&\tau,\quad &\mbox{in }\ \Omega,\\
u=&h,\quad &\mbox{on }\partial\Omega,
\ea
\right.
\eeq
with tension field $\tau\in \mathcal M^{1,a}(\Omega)$, for some $1< a<2$, and the boundary value
$h\in C^0(\partial\Omega, N)$,
and satisfies
\beq\label{vereg}
\int_{B_{2r_0}^+(x_0)}|\nabla u|^2\leq \ve_0^2,
\eeq
for some $x_0\in\partial\Omega$,  then $u\in C^0(B_{r_0}^+(x_0), N)$ and 
\beq\label{lcosc1}
 \underset{B_{r_0}^+(x_0)}{\rm{osc}}\ u
 \leq C\big(\int_{B_{2r_0}^+(x_0)}|\nabla u|^2\big)^{\frac12}
 +C\big(\|\tau\|_{\mathcal M^{1,a}(B_{2r_0}^+(x_0))}+\underset{T_{2r_0}(x_0)}{\rm{osc}} h\big).
\eeq
\end{theorem}

\pf With the help of Lemma \ref{interior1},  \eqref{lcosc1} can be proved similarly to that of 
Theorem 2.3. We only sketch it here.  By Fubini's theorem, there exists $r_1\in (\frac{3r_0}2,2r_0)$ such that 
\beq\notag
r_1\int_{S_{r_1}^+(x_0)}|\nabla u|^2\leq 8\int_{B_{2r_0}^+(x_0)}|\nabla u|^2.
\eeq
This, together with Sobolev's embedding theorem, implies that $u\in C^{\frac12}(S^+_{r_1}(x_0))$, and 
\beq\notag
\underset{S^+_{r_1}(x_0)}{\mbox{osc}}\ u
\leq C\big(r_1\int_{S^+_{r_1}(x_0)}|\nabla u|^2\big)^{\frac12}
\leq C\big(\int_{B^+_{2r_0}(x_0)}|\nabla u|^2\big)^{\frac12}.
\eeq 
Since $u=h$ on $T_{2r_0}(x_0)$, we then obtain that
\beq\label{lmeq1}
\underset{\partial B^+_{r_1}(x_0)}{\mbox{osc}}\ u
\leq \underset{T_{2r_0}(x_0)}{\mbox{osc}}h+C\big(\int_{B^+_{2r_0}(x_0)}|\nabla u|^2\big)^{\frac12}.
\eeq 
Now we  can apply Lemma \ref{interior1} and follow the same argument as in the proof of
(2.27) to show that 
\beq\label{lmeq2}
\underset{B_{r_0}^+(x_0)}{\mbox{osc}} u
\leq C\big(\underset{\partial B_{r_1}^+(x_0)}{\mbox{osc}} u
+\|\nabla u\|_{L^2(B_{2r_0}^+(x_0))}+\|\tau\|_{\mathcal M^{1,a}(B_{2r_0}^+(x_0))}\big).
\eeq
It is readily seen that \eqref{lcosc1} follows from \eqref{lmeq1} and \eqref{lmeq2}. The proof is now complete.
\endpf

\section{Proof of Theorem \ref{mainth1}}
\setcounter{equation}{0}
\setcounter{theorem}{0}
This section is devoted to the proof of Theorem \ref{mainth1}.  
In order to obtain both the energy identity \eqref{theq3} and the oscillation convergence \eqref{theq4}, we need to show,
similar to the proof of Theorem 1.2,  that there is neither energy concentration nor oscillation accumulation in the boundary neck regions.  More precisely, we need

\blm\label{lemma2}
For $\epsilon>0$, let $\{u_n\}\subset H^1(\Omega,N)$ be a sequence of approximate harmonic maps
with tension fields $\tau_n$ uniformly bounded in $L\log L(\Omega)\cap M^{1,a}(\Omega)$ for some
$1<a<2$ and Dirichlet boundary values $g_n$ satisfying (A1) and (A2).  
If, for sufficiently small $\delta>0$ and sufficiently large $R>1$, 
\beq\label{neckve}
\sup_{n\ge 1}\sup\limits_{Rr_n\leq \tau\leq 2\delta}\int_{B^+_{2\tau}(x_n)\setminus B^+_{\tau}(x_n)}|\nabla u_n|^2\leq \epsilon^2,
\eeq
holds for some $x_n\in\partial\Omega$ and $r_n\rightarrow 0$,
then there exists $\alpha=\alpha(a)>0$ such that
\beq\label{neeq1}
\int_{B_\delta^+(x_n)\setminus B_{4Rr_n}^+(x_n)}|\nabla u_n|^2\leq C(\epsilon^2+\delta+o(1)),
\eeq
where $\displaystyle\lim_{n\rightarrow\infty}o(1)=0$, and 
\beq\label{neeq2}
\underset{B_\delta^+(x_n)\setminus B_{2Rr_n}^+(x_n)}{\rm{osc}}\ u_n\leq C(\epsilon+\delta^\alpha).
\eeq
\elm

\pf After passing to a subsequence, we may assume that $x_n\rightarrow x_0\in\partial\Omega$.
For simplicity, we may assume that $x_n=x_0=0\in\partial\Omega$.  
For any $x\in B_\delta^+(0)\setminus B_{2Rr_n}^+(0)$, since
$$B_{\frac{|x|}2}(x)\cap\Omega\subset B_{\frac{3|x|}2}^+(0)\setminus B_{\frac{|x|}2}^+(0),$$
we have, by \eqref{neckve}, that 
\begin{equation}\label{neckve1}
\int_{B_{\frac{|x|}2}(x)\cap\Omega}|\nabla u_n|^2\le 2
\sup\limits_{Rr_n\leq \tau\leq 2\delta}\int_{B^+_{2\tau}(x_n)\setminus B^+_{\tau}(x_n)}|\nabla u_n|^2\le 2\epsilon^2.
\end{equation}
Observe that for any $a' \in (a, 2)$, it holds that $\mathcal M^{1,a}(B_{\frac{|x|}2}(x)\cap\Omega)
\subset \mathcal M^{1,a'}(B_{\frac{|x|}2}(x)\cap\Omega)$, and
$$\big\|f\big\|_{\mathcal M^{1, a'}(B_{\frac{|x|}2}(x)\cap\Omega)}
\le C|x|^{a'-a}\big\|f\big\|_{\mathcal M^{1, a}(B_{\frac{|x|}2}(x)\cap\Omega)},
\ \forall\ f\in \mathcal M^{1, a}(B_{\frac{|x|}2}(x)\cap\Omega).$$
Thus we can apply both Lemma \ref{interior1} and Theorem \ref{boundary1} to conclude that,
for some $a'\in (a,2)$,
\begin{eqnarray}\label{osc_estimate}
\underset{B_{\frac{|x|}4}(x)\cap\Omega} {\rm{osc}} u_n
&\le& C\big(\epsilon+\underset{B_{\frac{|x|}2}(x)\cap\partial\Omega}{\rm{osc}} h_n
+\|\tau_n\|_{\mathcal M^{1, a'}(B_{\frac{|x|}2}(x)\cap\Omega)}\big)\nonumber\\
&\le& C\big(\epsilon+\underset{B_{\frac{|x|}2}(x)\cap\partial\Omega}{\rm{osc}} h_n
+|x|^{a'-a}\|\tau_n\|_{\mathcal M^{1, a}(B_{\frac{|x|}2}(x)\cap\Omega)}\big)\nonumber\\
&\le& C\big(\epsilon+\underset{T_{2\delta}(0)}{\rm{osc}} h_n
+\delta^{a'-a}\big).
\end{eqnarray}
Note that the assumption (A1) implies that for sufficiently small $\delta>0$, it holds
$$\sup_{n\ge 1}\underset{T_{2\delta}(0)}{\rm{osc}} h_n\le \epsilon.$$
Therefore \eqref{osc_estimate} implies that
\beq\label{osc_estimate1}
\sup_{2Rr_n\le \tau\le \delta}\underset{S_\tau^+(0)}{\rm{osc}} u_n
\le C(\epsilon+\delta^{a'-a}).
\eeq
Since
$$\underset{B_\delta^+(0)\setminus B_{2Rr_n}^+(0)} {\rm{osc}} u_n
\le 2\sup_{2Rr_n\le \tau\le \delta}\underset{S_\tau^+(0)}{\rm{osc}} u_n
+\underset{T_{2\delta}(0)}{\rm{osc}} h_n,$$
we arrive at
\beq\label{osc_estimate2}
\underset{B_\delta^+(0)\setminus B_{2Rr_n}^+(0)} {\rm{osc}} u_n
\le C(\epsilon+\delta^{a'-a}).
\eeq
This yields \eqref{neeq2} with $\alpha=a'-a>0$.

Now we want to show \eqref{neeq1}.  From Fubini's theorem and the assumption \eqref{neckve}, 
we can choose $\delta_1\in (\frac{3\delta}4, \delta)$ so that
\beq\label{neeq3}
\delta_1\int_{S_{\delta_1}^+(0)}|\nabla u_n|^2
\leq 8\int_{B^+_{\delta}(0)\setminus B^+_{\frac{\delta}2}(0)}|\nabla u_n|^2\leq 8\epsilon^2.
\eeq
Let $G_n: B^+_{\delta_1}(0)\rightarrow \R^L$ be a solution of 
\beq\notag
\left\{
\ba{rll}
\Delta G_n=&0,\quad &\mbox{in } B^+_{\delta_1}(0),\\
G_n=&u_n,\quad &\mbox{on }\partial B^+_{\delta_1}(0).
\ea
\right.
\eeq
Then by the maximum principle and \eqref{osc_estimate2}, we have that
\beq\label{neeq4}
\underset{B^+_{\delta_1}(0)}{\mbox{osc}}\ G_n
\leq \underset{\partial B^+_{\delta_1}(0)}{\mbox{osc}}\ u_n
\leq C(\epsilon+\delta^\alpha).
\eeq
By the standard energy estimate of Laplace equation, \eqref{neeq3},  and assumption (A2) on $h_n$,
we also have
\beq\label{neeq5}
\begin{split}
&\int_{B^+_{\delta_1}(0)}|\nabla G_n|^2
\leq C\|u_n\|^2_{H^{\frac12}(\partial B_{\delta_1}^+(0))} \\
\leq &C\big(\|u_n\|^2_{H^{\frac12}(S_{\delta_1}^+(0))} 
+\|h_n\|^2_{H^{\frac12}(T_{\delta_1}(0))}\big)\\
\leq &C\big(\delta_1\|\nabla u_n\|^2_{L^2(S_{\delta_1}^+(0 ))} 
+\|h_n\|^2_{H^{\frac12}(T_{\delta}(0))}\big)\\
\leq & C\big(\|\nabla u_n\|^2_{L^2\big(B_{\delta}^+(0)\setminus B^+_{\frac{\delta}{2}}(0)\big)} 
+\|h_n\|^2_{H^{\frac12}(T_{\delta}(0))}\big)
\leq C(\epsilon+o(1)).
\end{split}
\eeq
We again apply Fubini's theroem to choose $\tau_n\in (2Rr_n, 4Rr_n)$ such that
\beq\label{neeq6}
\tau_n \int_{S_{\tau_n}^+(0)}|\nabla u_n|^2\le 4\int_{B_{4Rr_n}^+(0)\setminus B_{2Rr_n}^+(0)}|\nabla u_n|^2
\le 4\epsilon^2.
\eeq
Now we multiply the first equation of \eqref{nahm} by $u_n-G_n$, integrate over 
$B_{\delta_1}^+(0)\setminus B_{\tau_n}^+(0)$,  use the fact $u_n-G_n=0$ on $\partial B^+_{\delta_1}(0)$,
and apply \eqref{osc_estimate2}, \eqref{neeq3}, \eqref{neeq4},   \eqref{neeq6} to obtain
\beq\notag
\begin{split}
&\int_{B_{\delta_1}^+(0)\setminus B_{\tau_n}^+(0)}\nabla u_n\cdot\nabla (u_n-G_n)\\
=&\int_{\partial (B_{\delta_1}^+(0)\setminus B_{\tau_n}^+(0))}\frac{\partial u_n}{\partial\nu}(u_n-G_n)
+\int_{B_{\delta_1}^+(0)\setminus B_{\tau_n}^+(0)}\big(A(u_n)(\nabla u_n,\nabla u_n)-\tau_n\big)(u_n-G_n)\\
=&-\int_{S^+_{\tau_n}(0)}\frac{\partial u_n}{\partial r}(u_n-G_n)+\int_{B_{\delta_1}^+(0)\setminus B_{\tau_n}^+(0)}\big(A(u_n)(\nabla u_n,\nabla u_n)-\tau_n\big)(u_n-G_n)\\
\leq &\|u_n-G_n\|_{L^{\infty}(B_{\delta_1}^+(0)\setminus B_{\tau_n}^+(0))}\big(\tau_n\int_{S^+_{\tau_n}(0)}|\nabla u_n|^2\big)^{\frac12}\\
+&C\|u_n-G_n\|_{L^{\infty}(B^+_{\delta_1}(0)\setminus B^+_{\tau_n}(0))}\int_{B_{\delta_1}^+(0)\setminus B_{\tau_n}^+(0)}(|\nabla u_n|^2+|\tau_n|)
\\
\leq & C(\epsilon+\delta^\alpha).
\end{split}
\eeq
This, combined with \eqref{neeq5} and H\"older's inequality, implies that
\beq\notag
\begin{split}
&\int_{B_{\delta_1}^+(0)\setminus B_{\tau_n}^+(0)}|\nabla u_n|^2\\
\leq &\int_{B_{\delta_1}^+(0)\setminus B_{\tau_n}^+(0)}\nabla u_n\cdot\nabla G_n
+C(\epsilon+\delta^\alpha)\\
\leq & \frac12\int_{B_{\delta_1}^+(0)\setminus B_{\tau_n}^+(0)}|\nabla u_n|^2
+C\int_{B_{\delta_1}^+(0)\setminus B_{\tau_n}^+(0)}|\nabla G_n|^2+C(\epsilon+\delta^\alpha)\\
\leq &\frac12\int_{B_{\delta_1}^+(0)\setminus B_{\tau_n}^+(0)}|\nabla u_n|^2
+C(\epsilon+\delta^\alpha+o(1)).
\end{split}
\eeq
Hence we obtain that 
\beq\notag
\int_{B_{\delta_1}^+(0)\setminus B_{\tau_n}^+(0)}|\nabla u_n|^2\leq C(\epsilon+\delta^\alpha+o(1)),
\eeq
which clearly yields \eqref{neeq1}.  The proof is now complete.
\endpf

\medskip

Now we are ready to give a proof of Theorem \ref{mainth1}.

\medskip
\noindent{\bf Proof of Theorem} \ref{mainth1}:

The general scheme of proof is similar to that of Theorem \ref{bubbling1}, we will provide it for the completeness.
Without loss of generality, we can assume $1<a<2$. 
Note that $L\log L(\Omega)\cap M^{1,a}(\Omega)\hookrightarrow (L^1(\Omega), {\rm{weak\ }} L^1{\rm{-topology}})$  is
compact. Hence from the assumptions (1.13) and (A1) we can assume, after taking a subsequence, that
$$u_n\rightharpoonup u \ {\rm{in}}\ H^1(\Omega), \ \tau_n\rightharpoonup \tau \ {\rm{in}}\ L^1(\Omega),
\ h_n\rightarrow h \ {\rm{in}}\ C^0(\Omega). 
$$
As in the proof of Theorem \ref{bubbling1}, we can define the interior and boundary concentration sets by
$$\Sigma_1=\bigcap_{r>0}\big\{x\in\Omega: \liminf_{n\rightarrow\infty}\int_{B_r(x)}|\nabla u_n|^2\ge\epsilon_0^2\big\},
\ \Sigma_2=\bigcap_{r>0}\big\{x\in\partial\Omega: \liminf_{n\rightarrow\infty}\int_{B_r^+(x)}|\nabla u_n|^2\ge\epsilon_0^2\big\}.$$
It is well-known that for $i=1,2$, $\Sigma_i$ is a finite set of $m_i$ points, with
$m_i\le \frac{E_0}{\epsilon_0^2}$ and $E_0=\sup_{n\ge 1}\int_{\Omega}|\nabla u_n|^2$.
It follows from Lemma \ref{interior1} and Theorem \ref{boundary1} that 
$$u_n\rightarrow u \ {\rm{in}}\  H^1_{\rm{loc}}(\overline\Omega\setminus\Sigma)\cap C^0_{\rm{loc}}(\overline\Omega\setminus(\Sigma_1\cup\Sigma_2)),$$
and $u$ is an approximate harmonic map with tension field $\tau\in L\log L(\Omega)\cap M^{1,a}(\Omega)$
and Dirichlet boundary value $g\in H^\frac12(\partial\Omega, N)\cap C^0(\partial\Omega, N)$. Moreover,
by applying Lemma \ref{interior1} and Theorem \ref{boundary1} again, we see that $u\in C^0(\overline\Omega, N)$.

As in the proof of Theorem \ref{bubbling1}, we can find all possible bubbles generated by $\Sigma_1\subset\Omega_{\delta_1}$ for some $\delta_1>0$ to obtain nontrivial harmonic maps $\{\omega_i\}_{i=1}^{k_1}\subset C^\infty(\mathbb S^2, N)$, 
sequences of bubbling points 
$\{x_n^i\}_{i=1}^{k_1}\subset \Omega_{\delta_1}$ and scales $\{r_n^i\}$ with $r_n^i\rightarrow 0$ for $1\le i\le k_1$ such that
$$u_n(x_n^i+r_n^i\cdot)\rightarrow \omega_i \ {\rm{in}}\ H^1_{\rm{loc}}(\R^2)\cap C^0_{\rm{loc}}(\R^2), 
 \ 1\le i\le k_1.$$

For each point $y_0\in\Sigma_2$, we can find a sequence of points $\{y_n\}\subset\Omega\rightarrow y_0$ and $r_n\rightarrow 0$
such that
\beq\label{bdry_concen}
\int_{B_{r_n}(y_n)\cap\Omega}|\nabla u_n|^2
=\max\big\{\int_{B_{r_n}(x)\cap\Omega}|\nabla u_n|^2: x\in B_{\delta_1}^+(y_0)\big\}=\frac{\epsilon_0^2}{10^2}.
\eeq
Define the blow up sequence
$v_n(x)=u_n(y_n+r_n x):\Omega_n\equiv r_n^{-1}(B_{\delta_1}^+(y_0)\setminus\{y_n\}) \to N$.
Then $v_n$ is an approximate harmonic map with tension field
$\tau(v_n)(\cdot)=r_n^2\tau_n(y_n+r_n\cdot)$, satisfying $v_n(x)=g_n(y_n+r_n x)$ for $x\in \partial_0 \Omega_n\equiv r_n^{-1}(T_{\delta_1}(y_0)\setminus\{y_n\})$.

Let $z_n\in \partial\Omega$ be such that $d_n:=|z_n-y_n|={\rm{dist}}(y_n, \partial\Omega)$. Then we have

\smallskip
\noindent{\it Claim}.  $\lambda_n:=\frac{d_n}{r_n}\rightarrow \infty$ as $n\rightarrow\infty$. For, otherwise,
$\lim_{n\rightarrow\infty}\lambda_n=\lambda\in [0, \infty)$ so that $\Omega_n\rightarrow \R^2_{-\lambda}$ as $n\rightarrow\infty$. . It follows from \eqref{bdry_concen}, Lemma \ref{interior1}, and Theorem \ref{boundary1}
that there exists a nontrivial harmonic map $\omega_0\in H^1(\R^2_{-\lambda}, N)$ such that 
$$v_n\rightarrow \omega_0\ {\rm{in}}\ H^1_{\rm{loc}}(\R^2_{-\lambda})\cap C^0_{\rm{loc}}(\R^2_{-\lambda}).$$
It follows from the assumption (A1) that there exists a point $p_0\in N$ 
such that $g_n(y_n+r_n\cdot)\rightarrow p_0$ on $\partial\R^2_{-\lambda}$. 
Thus $\omega_0=p_0$ on $\partial\R^2_{-\lambda}$. However, it is well-known \cite{Lemaire} that any finite energy harmonic map $v: \R^2_{-\lambda}\to N$, with $v={\rm{constant}}$ on $\partial\R^2_{-\lambda}$, is
a constant map. We get a contradiction.
Hence the claim holds. 

It follows from the claim that $\Omega_n\rightarrow \R^2$ as $n\rightarrow \infty$, 
$\omega_0\in C^\infty(\mathbb S^2, N)$ is a nontrivial harmonic map, and  
$$v_n\rightarrow\omega_0 \ {\rm{in}}\ H^1_{\rm{loc}}(\R^2)\cap C^0_{\rm{loc}}(\R^2).$$
Repeating this procedure for finitely many times, we can find all possible boundary bubbles
$\{\omega_j\}_{j=k_1+1}^{k_2}\subset C^\infty(\mathbb S^2, N)$, and all corresponding 
blowing up points $\{x_n^j\}_{j=k_1+1}^{k_2}\subset\Omega\setminus \Omega_{\delta_1}$ and
scales $\{r_n^j\}_{j=k_1+1}^{k_2}$ with $r_n^j\rightarrow 0$  for $k_1+1\le j\le k_2$ such that
\eqref{theq2} holds for $1\le i<j\le k_2$. 

Concerning the bubbles generated by $\Sigma_1$, it follows from \cite{wang16} Theorem 1.2 that 
\beq\label{interior-energy_id30}
\lim_{n\rightarrow\infty}\int_{\Omega_{\delta_1}}|\nabla u_n|^2
=\int_{\Omega_{\delta_1}}|\nabla u|^2+\sum_{i=1}^{k_1}\int_{\mathbb S^2}|\nabla\omega_i|^2,
\eeq
and
\beq\label{interior-oscillation_conv30}
\lim_{n\rightarrow\infty}\big\|u_n-u-\sum_{i=1}^{k_1}\big(\omega_i(\frac{\cdot-x_n^i}{r_n^i})-\omega_i(\infty)\big)\big\|_{L^\infty(\Omega_{\delta_1})}=0.
\eeq
Thus, in order to prove \eqref{theq3} and \eqref{theq4}, it suffices to show that
\beq\label{interior-energy_id3}
\lim_{n\rightarrow\infty}\int_{\Omega\setminus\Omega_{\delta_1}}|\nabla u_n|^2
=\int_{\Omega\setminus\Omega_{\delta_1}}|\nabla u|^2+\sum_{j=k_1+1}^{k_2}\int_{\mathbb S^2}|\nabla\omega_j|^2,
\eeq
and
\beq\label{interior-oscillation_conv3}
\lim_{n\rightarrow\infty}\big\|u_n-u-\sum_{j=k_1+1}^{k_2}\big(\omega_j(\frac{\cdot-x_n^j}{r_n^j})-\omega_j(\infty)\big)\big\|_{L^\infty(\Omega\setminus\Omega_{\delta_1})}=0.
\eeq

As in the proof of Theorem \ref{bubbling1},  we can assume that $k=:k_2=k_1+1$, 
$\Sigma_2=\{x_{k}\}$, and $\omega_{k}$ is the only bubble at $x_{k}$.
Hence for any $\epsilon>0$, there exist sufficiently large $R>1$
and sufficiently small $0<\delta<\delta_1$ such that
\beq\label{necksmall2}
\sup_{Rr_n^{k}\le r\le \delta}\int_{(B_{2r}(x_n^{k})\setminus B_r(x_n^{k}))\cap\Omega}|\nabla u_n|^2\le\epsilon^2,
\ \forall\ n\ge 1.
\eeq
Let $\{y_n^{k}\}\subset\partial\Omega\rightarrow x_{k}$ be such that
$d_n^{k}=|x_n^{k}-y_n^{k}|={\rm{dist}}(x_n^{k}, \partial\Omega)\rightarrow 0$.
Then $\lambda_n^k=\frac{d_n^{k}}{r_n^{k}}\rightarrow \infty$.

It follows from the proof of \eqref{interior-energy_id3} in \cite{wang16} that
\beq\label{interior_energy_id6}
\lim_{R\rightarrow\infty}\lim_{n\rightarrow\infty}\int_{B_{d_n^{k}}(x_n^{k})
\setminus B_{Rr_n^{k}}(x_n^{k})}|\nabla u_n|^2=0,
\eeq
and
\beq\label{interior_oscillation_conv10}
\lim_{R\rightarrow\infty}\lim_{n\rightarrow\infty}\underset{B_{d_n^k}(x_n^k)\setminus B_{Rr_n^k}(x_n^k)}
{\rm{osc}} u_n=0.
\eeq
Observe that
\begin{eqnarray}\label{inclusion3}
\big(B_{\frac{\delta}4}(x_n^k)\cap\Omega\big)\setminus B_{d_n^k}(x_n^k)
= \big((B_{\frac{\delta}4}(x_n^k)\setminus B_{4d_n^k}(x_n^k))\cap\Omega\big)
\cup \big((B_{4d_n^k}(x_n^k)\cap\Omega)\setminus B_{d_n^k}(x_n^k)\big),
\end{eqnarray}
and
\begin{eqnarray}\label{inclusion4}
(B_{\frac{\delta}4}(x_n^k)\setminus B_{4d_n^k}(x_n^k))\cap\Omega
\subset B_{\frac{\delta}2}^+(y_n^k)\setminus B_{2d_n^k}^+(y_n^k)
\subset (B_{\delta}(x_n^k)\cap\Omega)\setminus B_{d_n^k}(x_n^k).
\end{eqnarray}
It follows from \eqref{necksmall2} that
\beq\label{neck_energy9}
\int_{(B_{4d_n^k}(x_n^k)\cap\Omega)\setminus B_{d_n^k}(x_n^k)}|\nabla u_n|^2\le 4\epsilon^2,
\eeq
and
\beq\label{neck_small3}
\sup_{d_n^k\le r\le {\delta}}\int_{B_{2r}^+(y_n^k)\setminus B_r^+(y_n^k)}|\nabla u_n|^2\le \epsilon^2.
\eeq
Thus we can apply Lemma \ref {lemma2} to obtain that
\beq\label{neck_energy10}
\int_{B^+_{\frac{\delta}2}(y_n^k)\setminus B_{2d_n^k}^+(y_n^k)}|\nabla u_n|^2\le C(\epsilon^2+\delta+o(1)),
\eeq
and
\beq\label{neck_osc7}
\underset{B_{\frac{\delta}2}(y_n^k)\setminus B_{2d_n^k}^+(y_n^k)}{\rm{osc}}u_n
\le C(\epsilon+\delta^\alpha),
\eeq
for some $\alpha\in (0,1)$.
It follows from \eqref{neck_energy9} and \eqref{neck_energy10}, and the inclusions \eqref{inclusion3} and
\eqref{inclusion4} that 
\beq\label{neck_energy11}
\int_{(B_{\frac{\delta}4}(x_n^k)\setminus B_{4d_n^k}(x_n^k))\cap\Omega}|\nabla u_n|^2\le C(\epsilon^2+\delta+o(1)),
\eeq
which, combined with \eqref{interior_energy_id1}, implies
\beq\label{neck_energy12}
\lim_{R\rightarrow\infty}\lim_{n\rightarrow\infty}\int_{(B_{\frac{\delta}4}(x_n^m)\cap\Omega)\setminus B_{Rr_n^m}(x_n^m)}|\nabla u_n|^2\le C(\epsilon^2+\delta).
\eeq
This yields \eqref{interior-energy_id3} after sending $\delta\rightarrow 0$, since $\epsilon>0$ is arbitrarily small.

From the assumption \eqref{neck_small3}, we can apply Lemma \ref{interior1} and Theorem \ref{boundary1}
to obtain that
\beq\label{neck_osc8}
\underset{B_{2d_n^k}^+(y_n^k)\setminus B_{d_n^k}^+(y_n^k)}{\rm{osc}}u_n\le C(\epsilon+\delta^\alpha).
\eeq
Combining \eqref{interior_oscillation_conv10} together with \eqref{neck_osc7} and \eqref{neck_osc8} implies
$$\lim_{R\rightarrow \infty}\lim_{n\rightarrow\infty}\underset{(B_{\frac{\delta}4}(x_n^k)\cap\Omega)\setminus B_{Rr_n^k}(x_n^k)}{\rm{osc}}u_n\le C(\epsilon+\delta^\alpha),$$
which yields \eqref{interior-oscillation_conv3} after sending $\delta\rightarrow 0$, since $\epsilon>0$ is arbitrarily small.
The proof of Theorem \ref{bubbling1} is now complete. \qed

\bigskip
\noindent{\bf Acknowledgement}. T. Huang is supported by NSF of Shanghai grant 16ZR1423800. 
C. Wang is partially supported by NSF 1522869.


\end{document}